\documentclass[journal]{IEEEtran}

\usepackage{epsfig}
\usepackage{amsfonts}
\usepackage{amsmath}
\usepackage{amssymb}
\usepackage{verbatim}
\usepackage{sgame}
\usepackage{color}
\usepackage{latexsym}

\newtheorem{assumption}{Assumption}
\newtheorem{remark}{Remark}
\newtheorem{algorithm}{Algorithm}
\newtheorem{definition}{Definition}
\newtheorem{theorem}{Theorem}
\newtheorem{lemma}{Lemma}

\begin{document}

\title{Decentralized Q-Learning \\ for Stochastic Teams and Games\thanks{Part of this paper is presented at the IEEE Conference on Decision and Control, December 2015. This research was partially supported by the Natural Sciences and Engineering Research Council of Canada (NSERC).
}}

\author{G\"{u}rdal~Arslan~and~Serdar~Y\"{u}ksel,~\IEEEmembership{Member,~IEEE}
\thanks{G. Arslan is with the Department of Electrical Engineering, University of Hawaii at Manoa, 440 Holmes Hall, 2540 Dole Street, Honolulu, HI 96822, USA. {\tt\small gurdal@hawaii.edu}}%
\thanks{S. Y\"{u}ksel is with the Department of Mathematics and Statistics, Queen's University, Kingston, Ontario, CANADA, K7L 3N6. {\tt\small   yuksel@mast.queensu.ca}}%
}

%

\maketitle

\begin{abstract}
There are only a few learning algorithms applicable
to stochastic dynamic teams and games which generalize Markov decision processes to decentralized stochastic control problems involving possibly self-interested decision makers. Learning in games is generally difficult because of the non-stationary environment in
which each decision maker aims to learn its optimal decisions
with minimal information in the presence of the other decision
makers who are also learning. In stochastic dynamic games, learning is
more challenging because, while learning, the decision makers
alter the state of the system and hence the future cost. In
this paper, we present decentralized Q-learning algorithms for
stochastic games, and study their convergence for the weakly
acyclic case which includes team problems as an important special case.
The algorithm is decentralized in that each decision maker has access to only its local information, the state information, and the local cost realizations; furthermore, it is completely oblivious to the presence of other decision makers.
We show that these algorithms converge to equilibrium policies almost surely in large classes of stochastic games.
\end{abstract}


\section{Introduction}

This paper aims at developing new learning algorithms with desirable
convergence properties for certain classes of stochastic games, which are discrete-time dynamic games in which the history can be summarized by a ``state'' \cite{fink1964equilibrium}.
More specifically, we focus on weakly acyclic stochastic games that can be used to
model cooperative systems. The chief merit of the paper lies in the fact
that learning takes place in stochastic games, which are truly dynamic games, as opposed to learning in repeated games in which the same
single-stage game is played in every stage.
In stochastic games, the policies selected by the decision makers not only impact
their immediate cost but also alter the stage-games to be played in the future
through the state dynamics.  Hence, our results are applicable to a
significantly broader set of applications.

The existing literature on learning in stochastic games is very
small in comparison with the literature on learning in repeated games. As the method of
reinforcement learning gained popularity in the
context of Markov decision problems, a surge of
interest in generalizing the method of reinforcement learning, in
particular
Q-learning algorithm \cite{q-learning}, to stochastic games has
led to a set of publications primarily in the computer science
literature; see \cite{sho} and the references therein.
In many of these publications, the authors
tend to assume that the real objective of the agents \footnote{The terms ``agent'' and ``decision maker'' are used interchangeably.}   is for some
reason to find and play an equilibrium strategy (and sometimes
this even requires agents to somehow agree on a particular
equilibrium strategy), and not necessarily to pursue their own
objectives.
Another serious issue is
that the multi-agent algorithms introduced in many of these recent
papers are not scalable since each agent needs to maintain estimates of its
Q-factors for each state/joint action pair and compute an
equilibrium at each step of the algorithm using the updated
estimates, assuming that the actions and objectives are
exchanged between all agents.

Standard Q-learning, which enables an agent to learn how to play
optimally in a single-agent environment, has also been applied to
very specific multi agent applications \cite{Tan,Sen}. Here, each
agent runs a standard Q-learning algorithm by ignoring the other
agents, and hence information exchange between agents and
computational burden on each agent are substantially lower than
aforementioned multi-agent extensions of Q-learning algorithm.
Also, standard Q-learning in a multi-agent environment makes sense
from individual bounded rationality point of view. However, no
analytical
results exist regarding the properties of standard Q-learning in a
stochastic game setting.

We should also mention several attempts to extend a well-known
learning algorithm called Fictitious Play (FP) \cite{bro,rob} to
stochastic games \cite{Claus,Schoenmakers,fp-stochastic-games}.
The joint action learning algorithm presented in \cite{Claus} would be
computationally prohibitive quickly as the number of
agents/states/actions grow. The algorithms presented in
\cite{Claus} are claimed to be convergent to an equilibrium in
single-state single-stage common interest games but without a
proof. The extension of FP considered in \cite{Schoenmakers}
requires each agent to calculate a stationary policy at each step
in response to the empirical frequencies of the stationary
policies calculated and announced by other agents in the past. The
main contribution of \cite{Schoenmakers} is to show that such FP
algorithm is not convergent even in the simplest 2x2x2 stochastic
game where there are two states and two agents with two moves for
each agent. The version of FP used in \cite{fp-stochastic-games}
is applicable only to zero-sum games (strictly adversarial
games).

Other related work includes
\cite{borkar,Bowling,Bowling-scalable}. In \cite{borkar}, a
multi-agent version of an actor-critic algorithm \cite{Konda} is
shown to be convergent to generalized equilibria in a weak sense
of convergence, whereas in \cite{Bowling} a policy iteration
algorithm is presented without rigorous results for stochastic
games. The algorithms given in \cite{borkar,Bowling} are rational
from individual agent perspective, however they require higher
level of data storing and processing than standard Q-learning. The
paper \cite{Bowling-scalable} uses the policy iteration algorithm
given in \cite{Bowling} in conjunction with certain approximation
methods to deal with a large state-space in a specific card-game
without rigorous results.

We should emphasize that our viewpoint is individual bounded
rationality and strategic decision making, that is, agents should
act to pursue their own objectives even in the short term using
localized information and reasonable algorithms. It is also
desired that agent strategies converge to an agreeable solution in
cooperative situations where agent objectives are aligned with
system designer's objective even though agents do not necessarily
strive for converging to a particular strategy.

The rest of the paper is organized as follows. In \S\ref{se:sg},
the model is introduced. In \S~\ref{se:lph}, the specifics
of the learning paradigm and the standard  Q-learning algorithm is discussed, followed by the presentation of our first Q-learning algorithm for stochastic games and its convergence properties.  Generalizations of our main
results in \S~\ref{se:lph} are presented in \S \ref{se:ext}.
This is followed by a simulation study in \S
\ref{simulationsection}. The paper is concluded with some final remarks
in \S \ref{conclusionsection}. Appendices contain the proofs of
the technical results in the paper.

\section{Stochastic Dynamic Games}
\label{se:sg}

\begin{figure}[ht]
\footnotesize
\centering
\includegraphics[scale=0.3]{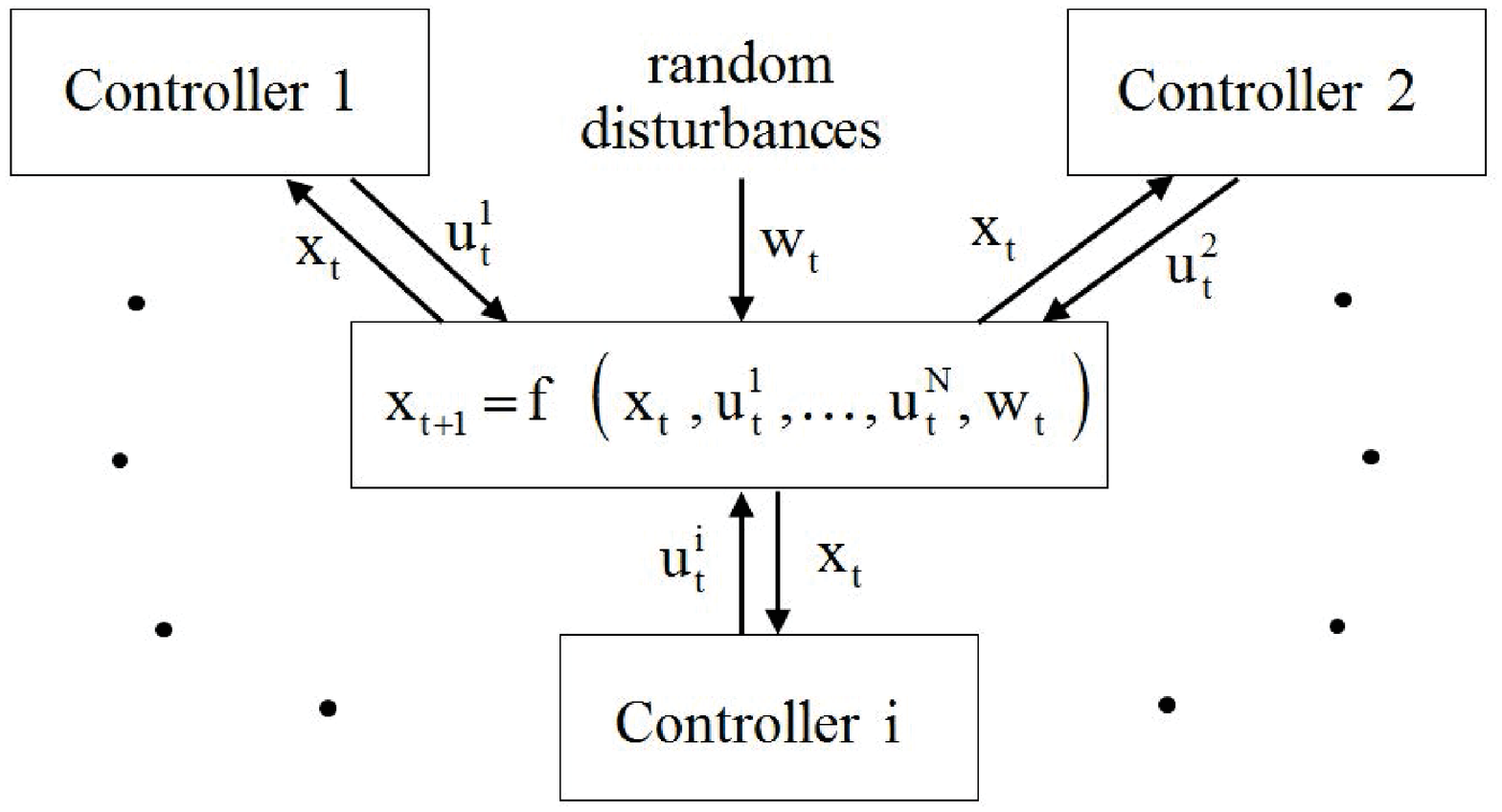}
\caption[]{A networked control system.}
\label{fig1}
\end{figure}

Consider the (discrete-time) networked control system illustrated in Figure~\ref{fig1}
where $x_t$ is the state of the system at time $t$, $u_t^i$ is the input generated by controller $i$ at time $t$,
and $w_t$ is the random disturbance input at time $t$.
Suppose that each controller $i$ is an autonomous decision maker (DM) interested in minimizing its own long-term cost
$$E \left[ \sum_{t\geq0 }  c^i(x_t,u_t^1,\dots,u_t^N) \right] $$
where $c^i(x_t,u_t^1,\dots,u_t^N)$ is the cost incurred by controller $i$ at time $t$, and $E[\cdot]$ denotes the expectation given a collection of control policies (which will be specified later in the paper) on a probability space $(\Omega,\mathcal{F},P)$. Although controller $i$ can only choose its own decisions $u_0^i,u_1^i,\dots$, its
cost generally depends on the decisions of all controllers through its single-stage cost as well as the state dynamics.
This dynamic coupling between self-interested DMs with long-term objectives naturally lead
to the framework of stochastic games \cite{fink1964equilibrium} which generalize Markov decision problems.

Over the past half-century, there have been many applications of stochastic games on control problems; see Chapter~XIV in \cite{bellmanadaptive} as an early reference. At the present time, the control theory literature includes a large number of papers employing the theory of stochastic games and their continuous-time counterparts called ``differential games'' \cite{isaacsdifferential}. Many papers in this body of work study a zero-sum game between a controller which aims to optimize the system performance and an adversary which controls certain system parameters and inputs to make the system performance as poor as possible. We selectively cite \cite{basbern} for robust control and minimax estimation problems, \cite{altman1994flow} for flow control in queueing networks, \cite{ding2013stochastic} for control of hybrid systems, and \cite{zhu2015game} for robustness, security, and resilience of cyber-physical control systems.
The case of nonzero-sum games in which the decision makers do not always have diametrically opposed objectives has also received significant attention; see for example \cite{altman1996non} on admission, service, and routing control in queueing systems, \cite{huang2010transmission} on transmission control in cognitive radio systems, \cite{zhu2010network} on network security, and \cite{gu2008differential} on formation control.

We should also mention the work on team decision problems where all DMs share a common long-term objective albeit with access to different information variables; see e.g., \cite{ho1980team,YukselBasarBook}.
In this paper, differently from the usual team decision problems in the literature, even though each DM has access to the state information, it does not have access to global information on the other DMs, and even their presence.
We also note that the emergence of distributed control systems requires the formulation of ``team problems'' within a game-theoretic framework where local controllers are tasked to achieve one system level objective without centralized coordination; see for example \cite{trodden2009distributed} on distributed model predictive control. This type of team problems and its generalizations where the objectives of DMs are aligned in some sense with a team objective are the primary focus of our work though the class of games considered in this paper is more general and it even includes some zero-sum stochastic games.

\subsection{Discounted Stochastic Dynamic Games}

A (finite) discounted stochastic game has the following ingredients; see \cite{fink1964equilibrium}.

\begin{itemize}
\item A finite set of DMs with the $i-$th DM referred to as DM$^i$ for $i\in\{1,\dots,N\}$
\item a finite set $\mathbb{X}$ of states
\item a finite set $\mathbb{U}^i$ of control decisions for each DM$^i$
\item a cost function $c^i$ for each DM$^i$ determining DM$^i$'s cost $c^i(x,u^1,\dots,u^N)$ at each state $x\in\mathbb{X}$ and for each joint decision $(u^1,\dots,u^N)\
    \in\mathbb{U}^1\times\cdots\mathbb{U}^N$
\item a discount factor $\beta^i\in(0,1)$ for each DM$^i$
\item a random initial state $x_0\in\mathbb{X}$
\item a transition kernel  for the probability
$P[ x^{\prime} |x,u^1,\dots,u^N]$ of each state transition from $x\in\mathbb{X}$ to $x^{\prime}\in\mathbb{X}$ for each joint decision $(u^1,\dots,u^N)\in\mathbb{U}^1\times\cdots\mathbb{U}^N.$

\end{itemize}

Such a stochastic game induces a discrete-time controlled
Markov process where the state at time $t$ is
denoted by $x_t\in\mathbb{X}$ starting with the initial state $x_0$.
At any time $t\geq0$, each DM$^i$ makes a control decision
$u_t^i\in \mathbb{U}^i$ (possibly randomly) based on the available information. The
state $x_t$ and the joint decisions $(u_t^1,...,u_t^N)$
together determine each DM$^i$'s cost $c^i(x_t,u_t^1,...,u_t^N)$ at time $t$ as well as the probability
distribution $P[ \ \cdot \ | \ x_t,u_t^1,...,u_t^N]$ with which the next state $x_{t+1}$ is
selected.

A policy for a DM is a rule of choosing an appropriate
control decision at any time based on the DM's history of observations. We will focus on stationary policies of the form where
a DM's decision at time $t$ is determined solely based on the state $x_t$.
Such policies for each DM$^i$ are identified by mappings from the state space $\mathbb{X}$ to the set $\mathcal{P}(\mathbb{U}^i)$ of probability distributions on $\mathbb{U}^i$. The interpretation is that a DM$^i$ using such a policy $\pi^i:\mathbb{X}\mapsto\mathcal{P}(\mathbb{U}^i)$ makes its decision $u_t^i$ at any time $t$
by choosing randomly from $\mathbb{U}^i$ according to $\pi^i(x_t)$. We will denote the set of such policies by $\Delta^i$ for each DM$^i$. We will primarily be interested in deterministic (stationary) policies\footnote{When it is not clear from the context, a ``policy'' will mean a deterministic policy.} denoted by $\Pi^i$ for each DM$^i$, where each  policy $\pi^i\in\Pi^i$ is identified by a mapping from $\mathbb{X}$ to $\mathbb{U}^i$.

The objective of each DM$^i$ is to find a policy
$\pi^i\in\Delta^i$ that minimizes its expected discounted  cost
\begin{equation}
J^i_{x}(\pi^1,\dots\pi^N)   =  E_{x} \left[   \sum_{t\geq0} (\beta^i)^t c^i\left(x_t,u^1_t,\dots,u^N_t\right) \right] \label{eq:dc}
\end{equation}
for all $x\in\mathbb{X}$, where $E_{x}$ denotes the conditional expectation given  $x_0=x$.
Since DMs have possibly different cost functions and each DM's
cost may depend on the control decisions of the other DMs, we adopt
the notion of equilibrium  to
represent those policies that are \textit{person-by-person optimal}. For ease of
notation, we denote the policies of all DMs other than DM$^i$ by $\pi^{-i}$. For future reference, we also define $\Pi^{-i}:=\times_{j\not=i} \Pi^j$ and $\Delta^{-i}:=\times_{j\not=i} \Delta^j$ as well as $\Pi:=\times_j \Pi^j$ and $\Delta:=\times_j \Delta^j$. Using this notation, we write a joint policy
$(\pi^1,\dots\pi^N)$ as $(\pi^{i},\pi^{-i})$ and $J_x^i(\pi^1,\dots\pi^N)$
as $J_x^i(\pi^{i},\pi^{-i})$.
\begin{definition}
A joint policy
$(\pi^{*1},\dots,\pi^{*N})\in\Delta$ constitutes
an (Markov perfect) equilibrium if, for all $i$, $x$,
$$J_x^i(\pi^{*i},\pi^{*-i})=\min_{\pi^i\in\Delta^i}J_x^i(\pi^{i},\pi^{*-i}).$$
\end{definition}
It is known that any finite discounted stochastic game possesses an equilibrium policy as defined above \cite{fud-tirole}.

Although the minimum above can always be achieved by a deterministic policy in $\Pi^i$ (since each DM$^i$'s problem is a stationary Markov decision problem when the policies of the other DMs are fixed at $\pi^{*-i}$), a deterministic equilibrium policy may not exist in general. However, many interesting classes of games do possess equilibrium in deterministic policies. In particular, large classes of games arising from applications where all DMs benefit from cooperation possess equilibrium in deterministic policies. The primary examples of such games of cooperation are team problems where all DMs have the same cost function. In team problems, the deterministic policies minimizing the common cost function are clearly equilibrium policies although non-optimal deterministic equilibrium policies may also exist. A more general set of games of cooperation are those in which some function, called the potential function, decreases whenever a single DM decreases its own cost by unilaterally switching from one deterministic policy to another one. In this class of games, the deterministic policies minimizing the potential function are equilibrium policies.  As such, we are primarily interested in the set of deterministic equilibrium policies denoted by $\Pi_{\rm eq}$, where $\Pi_{\rm eq}\subset\Pi$.

We next formally introduce the set of games considered in this paper.

\subsection{Weakly Acyclic Games}
\label{ss:wagbest}

Let $\Pi^i_{\pi^{-i}}$ denote DM$^i$'s set of (deterministic) best replies to any $\pi^{-i}\in\Delta^{-i}$, i.e.,
\begin{align*}
\Pi_{\pi^{-i}}^i  := \big\{\hat{\pi}^i\in\Pi^i  : \ & J_x(\hat{\pi}^i,\pi^{-i})=\min_{\pi^i\in\Delta^i}J_x(\pi^i,\pi^{-i}),  \\ & \ \mbox{for all} \  x \big\}.
\end{align*}
DM$^i$'s best replies to any $\pi^{-i}\in\Delta^{-i}$ can be characterized by its optimal Q-factors $Q_{\pi^{-i}}^i$ satisfying the fixed point equation
\begin{align}
\nonumber
Q_{\pi^{-i}}^i(x,u^i) = & E_{\pi^{-i}(x)} \big[c^i(x,u^i,u^{-i})  \\& +\beta^i \sum_{x^{\prime}\in\mathbb{X}} P[x^{\prime}|x,u^i,u^{-i}]\min_{v^i\in\mathbb{U}^i} Q_{\pi^{-i}}^i(x^{\prime},v^i) \big]
\label{eq:Qfp}
\end{align}
for all $x,u^i$, where $E_{\pi^{-i}(x)}$ denotes the expectation with respect to the joint distribution of $u^{-i}$ given by $\pi^{-i}(x)=\pi^1(x)\times\cdots\times\pi^{i-1}(x)\times\pi^{i+1}(x)\times\cdots\times\pi^N(x)$. The optimal Q-factor $Q_{\pi^{-i}}^i(x,u^i)$ represents DM$^i$'s expected discounted cost to go from the initial state $x$ assuming that DM$^i$ initially chooses $u^i$ and uses an optimal policy thereafter while the other DMs use $\pi^{-i}$.
One can then write  $\Pi_{\pi^{-i}}^i$ as
\begin{align*}
\Pi_{\pi^{-i}}^i = \big\{\hat{\pi}^i\in\Pi^i: & \ Q_{\pi^{-i}}^i(x,\hat{\pi}^i(x))=\min_{v^i\in\mathbb{U}^i} Q_{\pi^{-i}}^i(x,v^i),  \\ & \ \mbox{for all}  \  x \big\}.
\end{align*}
The set of (deterministic) joint best replies is denoted by $\Pi_{\pi}:=\Pi_{\pi^{-1}}^1\times\cdots\times\Pi_{\pi^{-N}}^N$.
Any best reply $\hat{\pi}^i\in\Pi_{\pi^{-i}}^i$ of DM$^i$ is called a \textit{strict best reply} with respect to $(\pi^i,\pi^{-i})$ if
\begin{align*}
J_x^i(\hat{\pi}^i,\pi^{-i}) & < J_x^i(\pi^i,\pi^{-i}), \quad\mbox{for some} \ x.
\end{align*}
Such a strict best reply $\hat{\pi}^i$ achieves DM$^i$'s minimum cost  given $\pi^{-i}$ for all initial states and results in a strict improvement over $\pi^i$ for at least one initial state.
\begin{definition}
We call a (possibly finite) sequence of deterministic joint policies $\pi_0,\pi_1,\dots$  a \textit{strict best reply path} if, for each $k$, $\pi_k$ and $\pi_{k+1}$ differ in exactly one DM position, say DM$^i$, and $\pi_{k+1}^i$ is a strict best reply with respect to $\pi_{k}$.
\end{definition}
\begin{definition}\label{def:best}
A discounted stochastic game is called \textit{weakly acyclic} under strict best replies if there is a strict best  reply path starting from each deterministic joint policy and ending at a deterministic equilibrium policy.
\end{definition}
\begin{figure}[ht]
\centering
\includegraphics[scale=0.23]{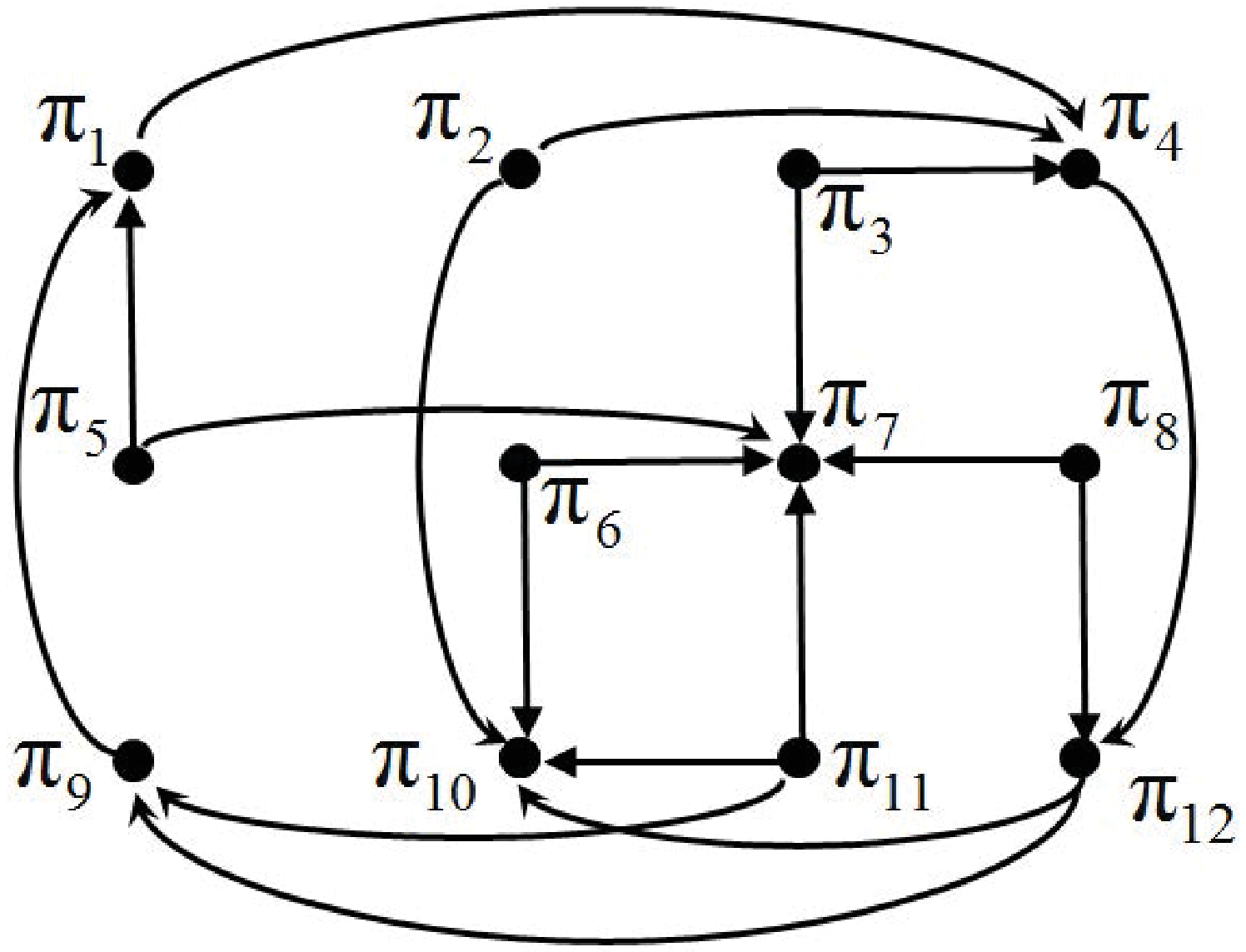}
\caption[]{The strict best reply graph of a stochastic game.} \label{fig3}
\end{figure}

Figure~\ref{fig3} shows the strict best reply graph of a game where the nodes represent the deterministic joint policies and the directed edges represent the single-DM strict best replies.
Each deterministic equilibrium policy is represented by a sink, i.e., a node with no outgoing edges, in such a graph. Note that the game illustrated in Figure~\ref{fig3} is weakly acyclic under strict best replies since there is a path from every node to a sink ($\pi_7$ or $\pi_{10}$). Note also that a weakly acyclic game may have cycles in its strict best reply graph, for example, $\pi_1\rightarrow\pi_4\rightarrow\pi_{12}\rightarrow\pi_9$ in Figure~\ref{fig3}.

Weakly acyclic games constitute a fairly large class of games. In the case of single-stage games, all potential games as well as dominance solvable games are examples of weakly acyclic games; see \cite{fabrikant2010structure}. We note that the concept of weak acyclicity introduced in this paper is with respect to the stationary Markov policies for stochastic games, and constitutes a generalization of weak acyclicity introduced in \cite{strategic-learning} for single-stage games. The primary examples of weakly acyclic games in the case of stochastic games are the team problems with finite state and control sets where DMs have identical cost functions and discount factors. Clearly, many other classes of stochastic games are weakly acyclic, e.g., appropriate multi-stage generalizations of potential games and dominance solvable games restricted to the stationary Markov policies are weakly acyclic for the same reason that the single-stage versions of these games are weakly acyclic \cite{fabrikant2010structure}.

\subsection{A Best Reply Process for Weakly Acyclic Games}
\label{ss:pap}
Consider a policy adjustment process in which only one DM updates its policy at each step by switching to one of its strict best replies. Such a process would terminate at an equilibrium policy if the game has no cycles in its strict best reply graph and the process continues until no DM has strict best replies. A weakly acyclic game may contain cycles in its strict best reply graph but there must be some edges leaving each cycle because otherwise there would not be a path from each node to a sink. Therefore, as long as each updating DM considers each of its strict best replies with positive probability, the adjustment process would terminate at an equilibrium policy in a weakly acyclic game with probability (w.p.) one. This adjustment process would require a criterion to determine the updating DM at each step and the DMs would have to a priori agree to this criterion. An equilibrium policy can be reached through a similar adjustment process without a pre-game agreement on the selection of the updating DM, if all DMs update their policies at each step but with some inertia.
Consider now the following policy adjustment process, which is the best reply process with memory length of one and inertia introduced in Sections 6.4-6.5 of  \cite{strategic-learning}.

\textit{Best Reply Process with Inertia (for DM$^i$):}

\noindent Set parameters\\
\hspace*{5mm} $\lambda^i\in(0,1)$: inertia\\
\noindent Initialize  $\pi_0^i \in \Pi^i$ (arbitrary) \\
Iterate $k\geq0$ \\
\hspace*{5mm} If $\pi_{k}^i\in\Pi_{\pi_k^{-i}}^i$ \\
\hspace*{10mm} $\pi_{k+1}^i=\pi_{k}^i$ \\
\hspace*{5mm} Else \\
\hspace*{10mm} $\pi_{k+1}^i = \left\{ \begin{array}{cl}   \pi_{k}^i & \mbox{ w.p. } \lambda^i \\ \mbox{any} \ \pi^i\in\Pi_{\pi_k^{-i}}^i & \mbox{ w.p. } (1-\lambda^i)/\Big|\Pi_{\pi_k^{-i}}^i\Big|\end{array}\right. $ \\
\hspace*{5mm} End \\

On the one hand, if the joint policy $\pi_k:=(\pi_k^1,\dots,\pi_k^N)$ is an equilibrium policy at any step $k$, then the policies will never change in the subsequent steps. On the other hand, regardless of what the joint policy $\pi_k:=(\pi_k^1,\dots,\pi_k^N)$ is at any step $k$, the joint policy $\pi_{k+L}$ in $L$ steps later will be an equilibrium policy with positive probability $p_{\min}>0$ where $L$ is the maximum length of the shortest strict best reply path from any policy to an equilibrium policy and  $p_{\min}$ depends only on the inertias $\lambda^1,\dots,\lambda^N$, and $L$. This readily implies that the best reply process with inertia will reach an equilibrium policy in finite number of steps w.p. $1$ \cite{strategic-learning}, i.e.,
$$P [\pi_k=\pi^*, \ \mbox{for some} \ \pi^*\in\Pi_{\rm eq} \ \mbox{and all large} \ k<\infty]=1.$$

We now note that each updating DM$^i$ at step $k$ needs to compute its best replies $\Pi_{\pi_k^{-i}}^i$, which can be done by first solving the fixed point equation (\ref{eq:Qfp}) for $\pi^{-i}=\pi_k^{-i}$. DM$^i$ can solve (\ref{eq:Qfp}), for example through value iterations, provided that DM$^i$ knows the state transition probabilities $P$ and the policies $\pi_k^{-i}$ of the other DMs  to evaluate the expectations in (\ref{eq:Qfp}). In most realistic situations, DMs would not have access to such information and therefore would not be able to compute their best replies directly. In the next section, we introduce our learning paradigm in which DMs would be able to learn their near best replies with minimal information and adjust their policies (approximately) along the strict best reply paths as in the best reply process with inertia.

\section{Q-Learning in Stochastic Dynamic Games}
\label{se:lph}
\subsection{Learning Paradigm for Stochastic Dynamic Games}
\label{se:lp}
The learning setup involves specifying the information that DMs have access to. We assume that each DM$^i$ knows its own set $\mathbb{U}^i$ of decisions and its own discount factor $\beta^i$. In addition, before choosing its decision $u_t^i$ at any time $t$, each DM$^i$ has the knowledge of
\begin{itemize}
\item its own past decisions $u_0^i,\dots,u_{t-1}^i$, and
\item past and current state realizations $x_0,\dots,x_t$, and
\item its own past cost realizations $$c^i(x_0,u_0^i,u_0^{-i}),\dots,c^i(x_{t-1},u_{t-1}^i,u_{t-1}^{-i}).$$
\end{itemize}
Each DM$^i$ has access to no other information such as the state transition probabilities or any information regarding the other DMs (not even the existence of the other DMs). In effect, the problem of decision making from the perspective of each DM$^i$ appears to be a stationary Markov decision problem. It is reasonable that each DM$^i$ with this view of its environment would use the standard Q-learning algorithm \cite{q-learning} to learn its optimal Q-factors and its optimal decisions. This would lead to the following Q-learning dynamics for each DM$^i$:
\begin{align*}
Q_{t+1}^i(x,u^i)  = & Q_t^i(x,u^i), \quad \mbox{for all} \ (x,u^i)\not= (x_t,u_t^i)\\
Q_{t+1}^i(x_t,u_t^i)  = & Q_t^i(x_t,u_t^i) + \alpha_t^i\big[c^i(x_t,u_t^i,u_t^{-i})
 \\& +\beta^i \min_{v^i\in\mathbb{U}^i} Q_t^i(x_{t+1},v^i)-Q_t^i(x_t,u_t^i)\big]
\end{align*}
where $\alpha_t^i\in[0,1]$ denotes DM$^i$'s step size at time $t$.

If only one DM, say DM$^i$, were to use Q-learning and the other DMs used constant policies $\pi^{-i}$, then DM$^i$ would asymptotically learn its corresponding optimal Q-factors, i.e., $$P [Q_t^i \rightarrow Q_{\pi^{-i}}^i]=1$$ provided that all state-control pairs $x,u^i$ are visited infinitely often and the step sizes are reduced at a proper rate. This follows from the well-known convergence of Q-learning in a stationary environment; see \cite{tsitsiklis1994asynchronous}.
To exploit the learnt Q-factors while maintaining exploration, the actual decisions are often selected with very high probability as
$$u_t^i \in \mbox{argmin}_{v^i\in\mathbb{U}^i} Q_t^i(x_t,v^i)$$
and with some small probability any decision in $\mathbb{U}^i$ is experimented. One common way of achieving this for DM$^i$ is to select any decision $u^i\in\mathbb{U}^i$ randomly according to (Boltzman action selection)
$$P [u_t^i = u^i|\mathcal{F}_t]= \frac{e^{- Q_t^i(x_t,u^i)/\tau}}{\sum_{v^i\in\mathbb{U}^i} e^{-Q_t^i(x_t,v^i)/\tau}} $$
where $\tau>0$ is a small constant called the temperature parameter, and $\mathcal{F}_t$ is the history of the random events realized up to the point just before the selection of $(u_t^1,\dots,u_t^N)$.

However, when all DMs use Q-learning and select their decisions as described above, the environment is non-stationary for all DMs, and there is no reason to expect convergence in that case. In fact, one can construct examples where DMs using Q-learning are caught up in persistent oscillations; see Section~4 in \cite{Leslie03b} for the non-convergence of Q-learning in Shapley's game. However, the convergence of Q-learning may still be possible in team problems, coordination-type games, or more generally in weakly-acyclic games. It is instructive to first consider the repeated games.

Here, there is no state dynamics (the set $\mathbb{X}$ of states is a singleton) and the DMs have no look-ahead ($\beta^1=\cdots\beta^N=0$). The only dynamics in this case is due to Q-learning which reduces to the averaging dynamics
\begin{align}
Q_{t+1}^i(u^i) & =  Q_t^i(u^i), \quad \mbox{for all} \ u^i\not= u_t^i \label{eq:qlrosg2}\\
Q_{t+1}^i(u_t^i) & =  Q_t^i(u_t^i) + \alpha_t^i\left[c^i(u_t^i,u_t^{-i})-Q_t^i(u_t^i)\right] \label{eq:qlrosg1}
\end{align}
where
\begin{equation}P [u_t^i = u^i|\mathcal{F}_t]= \frac{e^{- Q_t^i(u^i)/\tau}}{\sum_{v^i\in\mathbb{U}^i} e^{-Q_t^i(v^i)/\tau}}. \label{eq:qlrosg3} \end{equation}
The long-term behavior of these averaging dynamics is analyzed in \cite{Leslie03b} and strongly connected to the long-term behavior of the well-known Stochastic Fictitious Play (SFP) dynamics \cite{fudlev} in the case of two DMs; see Lemma~4.1 in \cite{Leslie03b}. In two-DM SFP, each DM$^i$ tracks the empirical frequencies of the past decisions of its opponent DM$^{-i}$ and chooses a nearly optimal response (with some experimentation) based on the incorrect assumption that DM$^{-i}$ will choose its decisions according to the empirical frequencies of its past decisions
$$
q_{t}^{-i}(u^{-i}) = \frac{1}{t} \sum_{k=0}^{t-1} I_{\{u_t^{-i}=u^{-i}\}}, \quad\mbox{for all} \ u^{-i} \\
$$
where $I_{\{\cdot\}}$ is the indicator function and
\begin{align*}
P [u_t^i = u^i|\mathcal{F}_t]  = & \frac{e^{- M_t^i(u^i)/\tau}}{\sum_{v^i\in\mathbb{U}^i} e^{-M_t^i(v^i)/\tau}} \\
M_t^i(u^i)  := & \sum_{u^{-i}} q_{t}^{-i}(u^{-i}) c^i(u^i,u^{-i}).
\end{align*}

Using the connection between Q-learning dynamics (\ref{eq:qlrosg1})-(\ref{eq:qlrosg3}) and SFP dynamics, the convergence of Q-learning (\ref{eq:qlrosg1})-(\ref{eq:qlrosg3}) is established in zero-sum games as well as in partnership games with two DMs; see Proposition~4.2 in \cite{Leslie03b}. It may be possible to extend this convergence result to multi-DM potential games \cite{potentialgames,Marden:TSMC:2009}, but this is currently unresolved. However, given the nonconvergence of FP (where DMs choose exact optimal responses with no experimentation, i.e., $\tau\downarrow0$) in some coordination games \cite{nonconv-coordination}, the prospect of establishing the convergence of Q-learning even in all two-DM weakly acyclic games does not seem promising.

It is possible to employ additional features such as the truncation of the observation history or multi-time-scale learning to obtain learning dynamics that are convergent in all repeated weakly acyclic games; see our own previous work \cite{Marden:JCO:2009} and the others \cite{Leslie03a,strategic-learning,Huck04,genweakenedfp}. However, the question of learning an equilibrium policy in stochastic games is an open question. The only relevant reference considering the stochastic games is \cite{borkar} where each DM uses value learning coupled with policy search at a slower time-scale. The results in \cite{borkar} apply to all stochastic games and therefore they are necessarily quite weak. Loosely speaking, the main result in \cite{borkar} shows that the limit points of certain empirical measures (weighted with the step sizes) in the policy space constitute ``generalized Nash equilibria'', which in particular does not imply convergence of learning to an equilibrium policy. In the next subsection, we propose a simple variation of Q-learning which converges to an equilibrium policy in all weakly acyclic stochastic games.

\subsection{Q-Learning in Stochastic Dynamic Games}
\label{ss:2tsq}

The discussion in the previous subsection reveals that the standard Q-learning (\ref{eq:qlrosg1})-(\ref{eq:qlrosg3}) can lead to robust oscillations even in repeated coordination games.
The main obstacle to convergence of Q-learning in games is due to the presence of multiple active learners leading to a non-stationary environment for all learners.  To overcome this obstacle, we use some inspiration from our previous work \cite{Marden:JCO:2009} on repeated games and modify the Q-learning for stochastic games as follows. In our variation of Q-learning, we allow DMs to use constant policies for extended periods of time called \textit{exploration phases}.

\begin{figure}[ht]
\footnotesize
\centering
\includegraphics[scale=0.23]{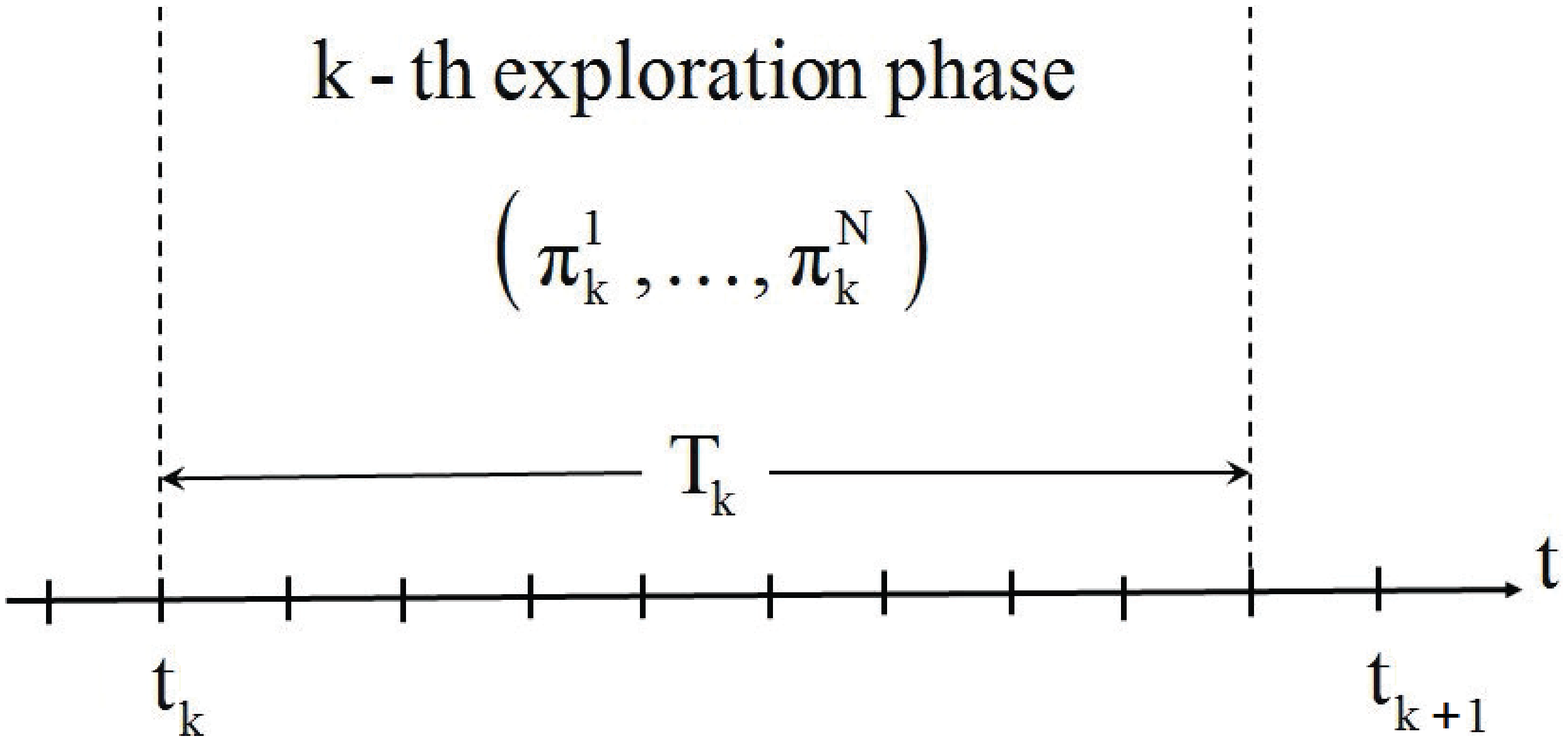}
\caption[]{An illustration of the $k-$th exploration phase.} \label{fig5}
\end{figure}
As illustrated in Figure~\ref{fig5}, the $k-$th exploration phase runs through times $t=t_k,\dots,t_{k+1}-1$, where $$t_{k+1}=t_k+T_k \qquad \mbox{(with $t_k=0$)}$$ for some integer $T_k\in[1,\infty)$ denoting the length of the $k-$th exploration phase. During the $k-$th exploration phase, DMs use some constant policies $\pi_k^1,\dots,\pi_k^N$ as their baseline policies with occasional experimentation.   The essence of the main idea is to create a stationary environment over each exploration phase so that DMs can accurately learn their optimal Q-factors corresponding to the constant policies used during each exploration phase.
Before arguing why this would lead to an equilibrium policy in all weakly acyclic stochastic games, let us introduce our variation of Q-learning more precisely.

\begin{algorithm}[for DM$^i$]
\label{al:1}

\noindent Set parameters\\
\hspace*{5mm} $\mathbb{Q}^i$: some  compact subset of the Euclidian space $\mathbb{R}^{|\mathbb{X}\times\mathbb{U}^i|}$ \\
\hspace*{11mm} where $|\mathbb{X}\times\mathbb{U}^i|$ is the number of pairs $(x,u^i)$\\
\hspace*{5mm} $\{T_k\}_{k\geq0}$: sequence of integers in $[1,\infty)$ \\
\hspace*{5mm} $\rho^i\in(0,1)$: experimentation probability\\
\hspace*{5mm} $\lambda^i\in(0,1)$: inertia\\
\hspace*{5mm} $\delta^i\in(0,\infty)$: tolerance level for sub-optimality\\
\hspace*{5mm} $\{\alpha_{n}^{i}\}_{n\geq0}$:  sequence  of step sizes where \\
\hspace*{20mm} $\alpha_{n}^{i}\in[0,1]$, $\sum_n\alpha_n^{i} = \infty$, $\sum_{n} \big(\alpha_n^{i}\big)^2 < \infty$ \\
\hspace*{20mm}    (e.g., $\alpha_n^i=1/n^r$ where $r\in(1/2,1]$)

\noindent Initialize  $\pi_0^i \in \Pi^i$ (arbitrary), $Q_0^i\in\mathbb{Q}^i$ (arbitrary) \\
Receive $x_0$ \\
Iterate $k\geq0$ \\
\hspace*{5mm}($k-$th exploration phase)\\
\hspace*{5mm}Iterate $t =t_k,\dots,t_{k+1}-1$ \\
\hspace*{10mm}  $u_t^i = \left\{\begin{array}{cl} \pi_k^i(x_t), & \textrm{ w.p. } 1-\rho^i\\ \textrm{any } u^i\in\mathbb{U}^i, & \textrm{ w.p. } \rho^i/|\mathbb{U}^i| \end{array} \right.$ \\
\hspace*{10mm}    Receive $c^i(x_t,u_t^i,u_t^{-i})$ \\
\hspace*{10mm}    Receive $x_{t+1}$ (selected according to  $P[ \ \cdot \ | \ x_t,u_t^i,u_t^{-i}]$) \\
\hspace*{10mm}    $n_t^i =$ the number of visits to $(x_t,u_t^i)$ in the $k-$th \\
\hspace*{19mm}    exploration phase up to $t$ \\
\hspace*{10mm} $Q_{t+1}^i(x_t,u_t^i)  =   (1-\alpha_{n_t^i}^{i})Q_t^i(x_t,u_t^i)$ \\
\hspace*{20mm} $ + \alpha_{n_t^i}^{i} \big[ c^i(x_t,u_t^i,u_t^{-i})  +  \ \beta^i \min_{v^i} Q_t^i(x_{t+1},v^i) \big]$\\
\hspace*{10mm}    $Q_{t+1}^i(x,u^i) =   Q_t^i(x,u^i)$,   for all $(x,u^i)\not=(x_t,u_t^i)$\\
\hspace*{5mm} End \\
\hspace*{5mm} $\Pi_{k+1}^i  = \big\{\hat{\pi}^i\in\Pi^i:  Q_{t_{k+1}}^i(x,\hat{\pi}^i(x))$ \\
\hspace*{35mm} $\leq \min_{v^i}Q_{t_{k+1}}^i(x,v^i)+\delta^i, \ \mbox{for all}  \ x\big\}$ \\
\hspace*{5mm} If  $\pi_k^i \in \Pi_{k+1}^i$ \\
\hspace*{10mm} $\pi_{k+1}^i = \pi_k^i$ \\
\hspace*{5mm} Else \\
\hspace*{10mm} $\pi_{k+1}^i = \left\{\begin{array}{cl} \pi_{k}^i, & \textrm{ w.p. } \lambda^i\\ \textrm{any } \pi^i\in\Pi_{k+1}^i, & \textrm{ w.p. } (1-\lambda^i)/|\Pi_{k+1}^i|\end{array} \right.  $  \\
\hspace*{5mm} End \\
\hspace*{5mm} Reset $Q_{t_{k+1}}^i$ to any $Q^i\in\mathbb{Q}^i$ (e.g., project $Q_{t_{k+1}}^i$ onto $\mathbb{Q}^i$)\\
End
\end{algorithm}

Algorithm~\ref{al:1} mimics the best reply process with inertia in \S\ref{ss:pap} arbitrarily closely with arbitrarily high probability under certain conditions. The key difference here is that each DM using Algorithm~\ref{al:1} approximately learns its optimal Q-factors during each exploration phase with limited observations. Accordingly, each DM updates its (baseline) policy to one of its near best replies with inertia based on its learnt Q-factors. Hence, Algorithm~\ref{al:1} can be regarded as an approximation to the best reply process with inertia in \S\ref{ss:pap}; see \cite{chapman2013convergent} where best replies are obtained based on rewards that must be estimated using noisy observations.


\begin{assumption}
\label{as:alpha}
For all $(x^{\prime},x)$, there exists a finite integer $H\geq0$ and joint decisions $\tilde{u}_0,\dots,\tilde{u}_{H}$ such that $$P[ x_{H+1}=x^{\prime} \ | \ (x_0,u_0,\dots,u_{H})=(x,\tilde{u}_0,\dots,\tilde{u}_{H})]>0.$$
\end{assumption}

Assumption~\ref{as:alpha}  ensures that the step sizes satisfy the well-known conditions of the stochastic approximation theory \cite{tsitsiklis1994asynchronous} during each exploration phase.

\begin{assumption}
\label{as:qp}
For all $i$, $0<\delta^i<\bar{\delta}$ and $0<\rho^i<\bar{\rho}$, where $\bar{\delta}$ and $\bar{\rho}$ (which depend only on the parameters of the game at hand) are defined in Appendix~\ref{Qsection}.
\end{assumption}

Assumption~\ref{as:qp} requires that the tolerance levels for sub-optimality used in the computation of near best replies as well as the experimentation probabilities be nonzero but sufficiently small.

\begin{theorem}\label{normalQ}
Consider a discounted stochastic game that is weakly acyclic under strict best replies. Suppose that  each DM$^i$ updates its policies by Algorithm~\ref{al:1}. Let Assumption~\ref{as:alpha}~and~\ref{as:qp} hold.
\begin{enumerate}
\item[(i)] For any $\epsilon>0$, there exist $\tilde{T}<\infty$, $\tilde{k}<\infty$ such that if $\min_{\ell}T_{\ell}\geq \tilde{T}$,  then
$$P\left[ \pi_{k} \in \Pi_{\rm eq} \right] \geq 1-\epsilon, \qquad \mbox{for all} \ k\geq\tilde{k}.$$
\item[(ii)]  If $T_k \rightarrow \infty$, then $$P\left[ \pi_{k} \in \Pi_{\rm eq} \right] \rightarrow 1.$$
\item[(iii)] There exists finite integers $\{\tilde{T}_k\}_{k\geq0}$ such that if $T_k\geq\tilde{T}_k$, for all $k$, then
$$ P\big[ \pi_{k} \rightarrow \pi^*, \ \mbox{for some} \ \pi^*\in\Pi_{\rm eq}\big] = 1.$$
\end{enumerate}
\end{theorem}

\begin{IEEEproof}
See Appendix~\ref{Qsection}.\qquad
\end{IEEEproof}

Let us discuss the main idea behind this result. Since all DMs use constant policies throughout any particular exploration phase, each DM indeed faces a stationary Markov decision problem in each exploration phase. Therefore, if the length of each exploration phase is long enough and the experimentation probabilities $\rho^1,\dots,\rho^N$ are small enough (but non-zero), each DM$^i$ can learn its corresponding optimal Q-factors in each exploration phase with arbitrary accuracy with arbitrarily high probability. This allows each DM$^i$ to accurately compute its near best replies to the other DMs' policies $\pi_k^{-i}$ at the end of the $k-$th exploration phase. Intuitively, allowing each DM$^i$ to update its policy $\pi_k^i$ to its near best replies (to $\pi_k^{-i}$) at the end of the $k-$th exploration phase with some inertia $\lambda^i\in(0,1)$ results in a policy adjustment process that approximates the best reply process with inertia in \S\ref{ss:pap}.

\begin{remark}
One may also wish to find explicit lower-bounds on $T_k$ to achieve almost sure convergence based on the convergence rates of the standard Q-learning with a single DM; we refer the reader to  \cite{even2004learning} for bounds on the convergence rates for standard Q-learning.
\end{remark}

\section{Generalizations}\label{se:ext}
\subsection{Learning in Weakly Acyclic Games under Strict Better Replies}
\label{ss:wagbetter}
We present another Q-learning algorithm with provable convergence to equilibrium in discounted stochastic games that are weakly acyclic under strict better replies. For this, we first introduce the notion of weak acyclicity under strict better replies. Given any $\pi=(\pi^i,\pi^{-i})\in\Delta$, let $\Upsilon^i_{\pi}$ denote DM$^i$'s set of (deterministic) better replies with respect to $\pi$, i.e.,
$$
\Upsilon_{\pi}^i := \big\{\hat{\pi}^i\in\Pi^i:J_x(\hat{\pi}^i,\pi^{-i}) \leq J_x(\pi^i,\pi^{-i}),  \ \mbox{for all} \  x \big\}.
$$
Any better reply $\hat{\pi}^i\in\Upsilon_{\pi}^i$ of DM$^i$ is called a \textit{strict better reply} (with respect to $\pi$) if
$$
J_x^i(\hat{\pi}^i,\pi^{-i})  < J_x^i(\pi^i,\pi^{-i}), \quad\mbox{for some} \ x.
$$
\begin{definition}
We call a (possibly finite) sequence of deterministic joint policies $\pi_0,\pi_1,\dots$  a \textit{strict better reply path} if, for each $k$, $\pi_k$ and $\pi_{k+1}$ differ in exactly one DM position, say DM$^i$, and $\pi_{k+1}^i$ is a strict better reply with respect to $\pi_{k}$.
\end{definition}
\begin{definition}\label{def:better}
A discounted stochastic game is called weakly acyclic under strict better replies if there is a strict better  reply path starting from each deterministic joint policy and ending at a deterministic equilibrium policy.
\end{definition}

Since every strict best reply path is also a strict better reply path, the set of games weakly acyclic under better replies contain (in fact, strictly) the set of games weakly acyclic under best replies.

It is straightforward to introduce a policy adjustment process analogous to the one in \S\ref{ss:pap} where, at each step, each DM$^i$ switches to one of its strict better replies with some inertia; see Sections 6.4-6.5 in \cite{strategic-learning}. Such a process would clearly converge to an equilibrium in games that are weakly acyclic under strict better replies. We next introduce a learning algorithm which allows each DM to learn the Q-factors corresponding to two policies, a baseline policy and a randomly selected experimental policy, during each exploration phase. If the learnt Q-factors indicate that the experimental policy is better than the baseline policy within a certain tolerance level, then the baseline policy is updated to the experimental policy with some inertia at the end of each exploration phase.
This learning algorithm enables DMs to adjust their policies with much less information (as in \S\ref{se:lp}), and follow (approximately) along the strict better reply paths that the adjustment process follows.

\begin{algorithm}[for DM$^i$]
\label{al:2}

\noindent Set parameters as in Algorithm~\ref{al:1}\\
\noindent Initialize  $\pi_0^i,\hat{\pi}_0^i \in \Pi^i$ (arbitrary except  $\hat{\pi}_0^i \not = \pi_0^i$), $Q_0^i,\hat{Q}_0^i\in\mathbb{Q}^i$ (arbitrary) \\
Receive $x_0$  \\
Iterate $k\geq0$ \\
\hspace*{5mm}($k-$th exploration phase)\\
\hspace*{5mm}Iterate $t =t_k,\dots,t_{k+1}-1$ \\
\hspace*{10mm}  $u_t^i = \left\{\begin{array}{cl} \pi_k^i(x_t), & \textrm{ w.p. } 1-\rho^i\\ \textrm{any } u^i\in\mathbb{U}^i, & \textrm{ w.p. } \rho^i/|\mathbb{U}^i|\end{array} \right.$ \\
\hspace*{10mm}    Receive $c^i(x_t,u_t^i,u_t^{-i})$ \\
\hspace*{10mm}    Receive $x_{t+1}$ (selected according to  $P[ \ \cdot \ | \ x_t,u_t^i,u_t^{-i}]$) \\
\hspace*{10mm}    $n_t^i =$ the number of visits to $(x_t,u_t^i)$ in the $k-$th \\
\hspace*{19mm}    exploration phase up to $t$ \\
\hspace*{10mm}    $Q_{t+1}^i(x_t,u_t^i)  =   (1-\alpha_{n_t^i}^{i})Q_t^i(x_t,u_t^i)$ \\
\hspace*{19mm}    $+ \alpha_{n_t^i}^{i} \big[ c^i(x_t,u_t^i,u_t^{-i})  +  \ \beta^i  Q_t^i(x_{t+1},\pi_k^i(x_{t+1})) \big]$\\
\hspace*{10mm}    $\hat{Q}_{t+1}^i(x_t,u_t^i)  =   (1-\alpha_{n_t^i}^{i})\hat{Q}_t^i(x_t,u_t^i) $\\
\hspace*{19mm}    $+ \alpha_{n_t^i}^{i} \big[ c^i(x_t,u_t^i,u_t^{-i})  +  \ \beta^i  Q_t^i(x_{t+1},\hat{\pi}_k^i(x_{t+1})) \big]$\\
\hspace*{10mm}    $Q_{t+1}^i(x,u^i) =   Q_t^i(x,u^i)$,   for all $(x,u^i)\not=(x_t,u_t^i)$\\
\hspace*{10mm}    $\hat{Q}_{t+1}^i(x,u^i) =  \hat{Q}_t^i(x,u^i)$,   for all $(x,u^i)\not=(x_t,u_t^i)$\\
\hspace*{5mm} End \\
 \hspace*{5mm} If  $(\hat{Q}_{t_{k+1}}^i(x,\hat{\pi}_k^i(x))  \leq Q_{t_{k+1}}^i(x,\pi_k^i(x)) + \delta^i, \ \mbox{for all} \ x)$\\
\hspace*{45mm} and \\
\hspace*{5mm}     \ \ \ $(\hat{Q}_{t_{k+1}}^i(x,\hat{\pi}_k^i(x))  \leq Q_{t_{k+1}}^i(x,\pi_k^i(x)) - \delta^i, \ \mbox{for some} \ x)$
 \\
\hspace*{10mm} $\pi_{k+1}^i = \left\{\begin{array}{cl} \pi_{k}^i, & \textrm{ w.p. } \lambda^i\\  \hat{\pi}_k^i, & \textrm{ w.p. } 1-\lambda^i\end{array} \right.$  \\
\hspace*{5mm} Else \\
\hspace*{10mm} $\pi_{k+1}^i = \pi_k^i$ \\
\hspace*{5mm} End \\
\hspace*{5mm} $\hat{\pi}_{k+1}^i=$ any policy  $\pi^i\in\Pi^i \backslash \{\pi_{k+1}^i\}$ with equal\\
\hspace*{17mm} probability\\
\hspace*{5mm} Reset $Q_{t_{k+1}}^i$, $\hat{Q}_{t_{k+1}}^i$ to any $Q^i,\hat{Q}^i\in\mathbb{Q}^i$\\
End
\end{algorithm}

Since any policy except the baseline policy can be chosen as an experimental policy (with equal probability),  each DM can switch to any of its strict better replies with positive probability. In contrast, each DM using Algorithm~\ref{al:1}  can only switch to one of its strict best replies. As a result, each DM using Algorithm~\ref{al:2} can escape a strict best reply cycle by switching to a strict better reply (if one exists); whereas, any DM using Algorithm~\ref{al:1} cannot. This flexibility comes at the cost of running two Q-learning recursions, one for the baseline policy and the other for the experimental policy, instead of one. However, this flexibility also leads to convergent behavior in a strictly larger set of games. We cite \cite{arslan2007equilibrium} as a reference to an earlier use of the idea of comparing two strategies and selecting one according to the Boltzman distribution.

The counterpart of Theorem~\ref{normalQ} can be obtained for Algorithm~\ref{al:2} in games that are weakly acyclic under strict better replies.

\begin{assumption}
\label{as:qp2}
For all $i$, $0<\delta^i<\check{\delta}$ and $0<\rho^i<\check{\rho}$, where $\check{\delta}$ and $\check{\rho}$ (which depend only on the parameters of the game at hand) are defined in Appendix~\ref{Qsection2}.
\end{assumption}

\begin{theorem}\label{normalQ2}
Consider a discounted stochastic game that is weakly acyclic under strict better replies. Suppose that each DM$^i$ updates its policies by Algorithm~\ref{al:2}. Let Assumption~\ref{as:alpha}~and~\ref{as:qp2} hold.
\begin{enumerate}
\item[(i)]
For any $\epsilon>0$, there exist $\tilde{T}<\infty$, $\tilde{k}<\infty$ such that if $\min_{\ell}T_{\ell}\geq \tilde{T}$, then
$$P\left[ \pi_{k} \in \Pi_{\rm eq} \right] \geq 1-\epsilon, \qquad  k\geq\tilde{k}.$$
\item[(ii)] If  $T_k \rightarrow \infty$,  then
$$P\left[ \pi_{k} \in \Pi_{\rm eq} \right] \rightarrow  1.$$
\item[(iii)] There exists finite integers $\{\tilde{T}_k\}_{k\geq0}$ such that if $T_k\geq\tilde{T}_k$, for all $k$, then
$$ P\big[ \pi_{k} \rightarrow  \pi^*, \ \mbox{for some} \ \pi^*\in\Pi_{\rm eq}  \big] = 1.$$
\end{enumerate}
\end{theorem}

\begin{IEEEproof}
See Appendix~\ref{Qsection2}.\qquad
\end{IEEEproof}

\subsection{Learning in Weakly Acyclic Games under multi-DM Strict Best or Better Replies}

The notion of weak acyclicity can be generalized by allowing multiple DMs to simultaneously update their policies in a strict best or better reply path.

\begin{definition}
We call a (possibly finite) sequence of deterministic joint policies $\pi_0,\pi_1,\dots$  a \textit{multi-DM strict best (better) reply path} if, for each $k$, $\pi_k$ and $\pi_{k+1}$ differ for at least one DM and, for each deviating DM$^i$,  $\pi_{k+1}^i$ is a strict best (better) reply with respect to $\pi_{k}$.
\end{definition}

\begin{definition}
A discounted stochastic game is called weakly acyclic under multi-DM strict best (better) replies if there is a multi-DM strict best (better) reply path starting from each deterministic joint policy and ending at a deterministic equilibrium policy.
\end{definition}

This generalization leads to a strictly larger set of games that are weakly acyclic. To see this, consider a single-stage game characterized by the cost matrices in Figure~\ref{fig:ewa} where DM$^1$ chooses a row, DM$^2$ chooses a column, and DM$^3$ chooses a matrix, simultaneously.
\begin{figure}[ht]
\footnotesize
\centering
\begin{game}{3}{3}[$1$]
& $1$ & $2$ & $3$\\
$1$ &$-a,0,0$ &$0,a,0$ &$0,-a,-a$\\
$2$ &$a,0,0$ &$-a,-a,0$ &$a,0,0$ \\
$3$ &$0,-a,-a$ & $0,a,0$ &$-a,0,-a$
\end{game}
\quad
\begin{game}{3}{3}[$2$]
& $1$ & $2$ & $3$\\
$1$ &$0,-a,-a$ &$0,0,0$ &$0,0,0$\\
$2$ &$a,0,0$ &$-a,0,-a$ &$-a,-a,-a$ \\
$3$ &$-a,-a,0$ & $0,0,0$ &$0,0,0$
\end{game}
\caption[]{Cost matrices of a single-stage game with three DMs.}\label{fig:ewa}
\end{figure}
Assume $a>0$. There is no strict best (or better) reply path to an equilibrium from the joint decisions $(1,1,1)$, $(1,3,1)$, $(3,3,1)$, $(3,1,1)$, $(1,1,2)$, $(3,1,2)$, if only a single DM can update its decision at a time. Therefore, this game is not weakly acyclic under strict best (or better) replies in the sense of Definition~\ref{def:best} (or Definition~\ref{def:better}). However, if multiple DMs are allowed to switch to their strict best (or better) replies simultaneously, then it becomes possible to reach the equilibrium $(2,3,2)$ from any joint decision. For example, if DM$^2$ and DM$^3$ switch to their strict best (or better) replies simultaneously from the joint decision $(1,1,1)$, then the resulting joint decision would be $(1,3,2)$. This would subsequently lead to the equilibrium $(2,3,2)$ if DM$^1$ switches to its strict best (or better) reply from $(1,3,2)$.

All learning algorithms introduced in the paper allow multiple DMs to simultaneously update their policies with positive probability. In view of this, it is straightforward to see that our main convergence results Theorem~\ref{normalQ} (Theorem~\ref{normalQ2}) hold in games that are weakly acyclic under multi-DM strict best (better) replies.

\section{A Simulation Study: Prisoner's Dilemma with a State}\label{simulationsection}
We consider a discounted stochastic game with two DMs where $\mathbb{X}=\mathbb{U}^1=\mathbb{U}^2=\{1,2\}$. Each DM$^i$'s utility (to be maximized) at each time $t\geq0$ depends only on the joint decisions $(u_t^1,u_t^2)$ of both DMs as
\begin{figure}[htb]
\footnotesize
\centering
\begin{game}{3}{2}[DM$^i$:][DM$^{-i}$:]
              & 1 & 2 \\
1     &$c$            & $a$ \\
2     &$b$            & $0$
\end{game}
\caption[]{DM$^i$'s single-stage utility.}
\end{figure}

\noindent We assume $b>c>0>a$. The state evolves as
\begin{align*}
P\big [x_{t+1}= 1 \ |  \ (u_t^1,u_t^2)=(1,1)\big]  & = 1-\gamma \\
P\big [x_{t+1}= 2 \ |  \ (u_t^1,u_t^2)\not=(1,1)\big] & = 1-\gamma
\end{align*}
where $\gamma\in(0,1)$ and  $P[x_0=1]=1/2$.

The single-stage game corresponds to the well-known prisoner's dilemma where the $i-$th prisoner (DM$^i$) cooperates (defects) at time $t$ by choosing $u_t^i=1$ ($u_t^i=2$). The single-stage game has a unique equilibrium $(u^1,u^2)=(2,2)$, i.e., both DMs defect, leading to rewards $(0,0)$. The dilemma is that each DM can do strictly better by cooperating, i.e., $(u^1,u^2)=(1,1)$ (not an equilibrium).

In the multi-stage game, the state $x_t$ indicates, w.p. $1-\gamma$, whether or not both DMs cooperated in the previous stage. It turns out that cooperation can be obtained as an equilibrium of the multi-stage game if the DMs are patient, i.e., the discount factors are sufficiently high, and the error probability $\gamma$ is sufficiently small . Note that each DM$^i$ has four different policies of the form $\pi^i:\mathbb{X}\rightarrow\mathbb{U}^i$. For large enough $\beta^1,\beta^2$, and small enough $\gamma$, the multi-stage game has two (Markov perfect) equilibria. In one equilibrium, called the cooperation equilibrium, each DM cooperates if $x=1$ and defects otherwise. In the other equilibrium, called the defection equilibrium, both DMs always defect. Furthermore, from any joint policy in $\Pi^1\times\Pi^2$, there is a strict best reply path to one of these two equilibria, which implies that the multi-stage game is weakly acyclic under strict best replies.

We set $b=2$, $c=1$, $a=-1$, $\gamma=0.3$. We simulate Algorithm~\ref{al:1} with the following parameter values: $\rho^i=0.1$, $\lambda^i=0.5$, $\delta^i=0$, $\alpha_k^i = 1/k^{0.51}$, for all $i,k$.
We keep the lengths of the exploration phases constants, i.e, $T_k=T$, for all $k$. We consider different values for $T$ since the lengths of the exploration phases appear to be most critical for the behavior of the learning process. For each value of $T$, we run Algorithm~\ref{al:1} and the best reply process with inertia (in \S\ref{ss:pap}) in parallel, with $1000$ policy updates starting from each of the $16$ initial joint policies in $\Pi$.
We initialize all the learnt Q-factors at $0$ for each simulation run; however, we do not reset the learnt Q-factors at the end of any exploration phase during any simulation run. We let $\pi_k$ and $\breve{\pi}_k$ denote the policies generated by Algorithm~\ref{al:1} and the best reply process with inertia in \S\ref{ss:pap}, respectively. For each value of $T$, Table~\ref{tb:1} shows the fraction of times at which $\pi_k$ visits an equilibrium and the fraction of times at which $\pi_k$ agrees with $\breve{\pi}_k$, during the $1000$ policy updates (averaged uniformly over all $16$ initial policies in $\Pi$).

The results in Table~\ref{tb:1} reveals that, as $T$ increases, $\pi_k$ visits an equilibrium and agrees with $\breve{\pi}_k$ more often. This is consistent with Theorem~\ref{normalQ} since DMs are expected to learn their Q-factors more accurately with higher probability for larger values of $T$. When $T$ is sufficiently large, the polices $\pi_k$ are at equilibrium and agrees with $\breve{\pi}_k$ nearly all of the time regardless of the initial policy. In a typical simulation run (with a large enough $T$), the polices $\pi_k$ and $\breve{\pi}_k$ transition to an equilibrium in few steps and stay at equilibrium thereafter.
\begin{table}[ht]
\begin{center}
    \begin{tabular}{ | c | c | c |}
    \hline
    $T$  & \begin{tabular}{c} $\frac{1}{1001}\sum_{k=0}^{1000} I_{\{\pi_k\in\Pi_e\}}$ \\ (averaged over $\pi_0\in\Pi$) \end{tabular} &
               \begin{tabular}{c} $\frac{1}{1001}\sum_{k=0}^{1000} I_{\{\pi_k=\breve{\pi}_k\}} $ \\ (averaged over $\pi_0\in\Pi$) \end{tabular} \\ \hline
    $10$ & $0.2581$ & $0.1254$\\ \hline
    $25$   & $0.5274$ & $0.3410$\\ \hline
    $50$   & $0.7835$ & $0.6170$\\ \hline
    $100$ & $0.9282$ & $0.6301$\\ \hline
    $1000$ & $0.9935$ & $0.6879$\\ \hline
    $10000$ & $0.9978$ & $0.7733$\\ \hline
    $50000$ & $0.9976$ & $0.9705$\\ \hline
    \end{tabular}
\end{center}
\label{tb:1}
\caption[]{The fraction of times at which $\pi_k$ visits an equilibrium and \\ the fraction of times at which $\pi_k$ agrees with $\breve{\pi}_k$.}
\end{table}

\section{Concluding Remarks}\label{conclusionsection}

In this paper, we develop  decentralized Q-learning algorithms and present their convergence properties for stochastic games under weak acyclicity. This is the first paper, to our knowledge, that presents learning algorithms with convergence to equilibria in large classes of stochastic games.  The decision makers observe only their own decisions and cost realizations, and the state transitions; they need not even know the presence of the other decision makers.

Our approach has a two-time scale flavor; however, unlike the existing work on multi-time-scale learning, it does not depend on the stochastic approximation theory. Note that the existing work on multi-time-scale learning, e.g., \cite{borkar,Konda,Leslie03b,Leslie03a},  require the stability analysis of some ordinary differential equations (ODE) describing the mean behavior of the learning algorithms. Aside from the difficulty of choosing the step sizes running at multiple time scales, the existing work involves nonlinear ODEs whose analysis does not seem to be within reach even for dynamic team problems. In contrast, our approach  leads to a considerably simpler analysis for all weakly acyclic stochastic games.

\appendices

\section{A Uniform Convergence Result for the Standard Q-Learning Algorithm with a Single DM}
\label{se:ap3}

Convergence of the standard Q-learning algorithm with a single DM is well known \cite{tsitsiklis1994asynchronous}. However, to prove the results of this paper, we need the sample paths generated by the standard Q-learning algorithm to well behave with respect to the initial conditions. Let us now consider a single-DM version of the setup introduced in \S\ref{se:sg} where the DM index $i$ (in the superscript) is dropped (only in Appendix~\ref{se:ap3}) and $c(x,u)$ representing the one-stage cost for applying control $u$ at $x$ is an exogenous random variable with finite variance. Let us assume that a single DM using a stationary random policy $\pi\in\Delta$ updates its Q-factors as: for $t\geq0$,
    \begin{align}
     Q_{t+1}(x,u)  & =   Q_t(x,u), \qquad\mbox{for all} \ (x,u)\not=(x_t,u_t) \label{eq:sQ1}\\
     Q_{t+1}(x_t,u_t)  & =   Q_t(x_t,u_t)    +  \alpha_{n_t} \Big( c(x_t,u_t) \nonumber \\ & \qquad +\beta \min_{v} Q_t(x_{t+1},v) -Q_t(x_t,u_t)   \Big) \label{eq:sQ2}
    \end{align}
where the initial condition $Q_0$ is given, $u_t$ is chosen according to $\pi(x_t)$, the state $x_t$ evolves according to $P[ \ \cdot \ | x_t,u_t]$ starting at $x_0$, $n_t$ is the number of visits to $(x_t,u_t)$ up to time $t$, and $\{\alpha_n\}_{n\geq0}$ is a sequence  of step sizes satisfying
    $$\alpha_n\in[0,1], \quad \sum_{n} \alpha_n = \infty, \quad  \quad \sum_{n} \alpha_n^2 < \infty.$$

\begin{lemma}
\label{lm:ql}
Assume that each $(x,u)$ is visited infinitely often w.p. $1$. For any $\epsilon>0$ and compact $\mathbb{Q}\in\mathbb{R}^{|\mathbb{X}\times\mathbb{U}|}$, there exists $T_{\epsilon}^{\mathbb{Q}}<\infty$ such that, for any $Q_0\in\mathbb{Q}$,
$$P\left[ \sup_{t\geq T_{\epsilon}^{\mathbb{Q}}} \left|Q_t - \bar{Q}\right|_{\infty} \leq \epsilon  \right] \geq 1-\epsilon$$
where $|\cdot|_{\infty}$ denotes the maximum norm and $\bar{Q}$ is the unique fixed point of the mapping $F:\mathbb{X}\times\mathbb{U}\mapsto\mathbb{X}\times\mathbb{U}$ defined by
$$F(Q)(x,u) = E[c(x,u)]+\beta\sum_{x^{\prime}} P[x^{\prime}|x,u] \min_v Q(x^{\prime},v)$$
for all $x$, $u$.
\end{lemma}

\begin{IEEEproof}
Let $\{Q_t^{\prime}\}_{t\geq0}$ and $\{Q_t^{\prime\prime}\}_{t\geq0}$ be the trajectories for the initial conditions $Q_0^{\prime}$ and $Q_0^{\prime\prime}$, respectively, corresponding to a sample path $\{(x_t,u_t,c(x_t,u_t))\}_{t\geq0}$. It is easy to see that, for all  $t\geq0$,
        \begin{align*}
     & |Q_{t+1}^{\prime}(x_t,u_t)-Q_{t+1}^{\prime\prime}(x_t,u_t)|  \\ & \ \ \leq   (1-\alpha_{n_t})|Q_{t}^{\prime}(x_t,u_t)-Q_{t}^{\prime\prime}(x_t,u_t)|    + \alpha_{n_t}\beta|Q_t^{\prime}-Q_t^{\prime\prime}|_{\infty}.
    \end{align*}
This implies that $M_t:=\sup_{Q_0^{\prime},Q_0^{\prime\prime}\in\mathbb{Q}}|Q_t^{\prime}-Q_t^{\prime\prime}|_{\infty}$  is non-increasing and therefore convergent.
Suppose that $M_t\rightarrow \bar{M}>0$. There exists some $\bar{t}<\infty$ such that $\max_{t\geq\bar{t}}M_t<\bar{M}(1+1/\beta)/2$. Hence, we have, for all $t\geq\bar{t}$,
        \begin{align*}
     & |Q_{t+1}^{\prime}(x_t,u_t)-Q_{t+1}^{\prime\prime}(x_t,u_t)| \\ & \  \leq   (1-\alpha_{n_t})|Q_{t}^{\prime}(x_t,u_t)-Q_{t}^{\prime\prime}(x_t,u_t)|  +\alpha_{n_t} \beta \frac{\bar{M}(1+1/\beta)}{2}.
    \end{align*}
This leads to: for all $(x,u)$ and $t\geq\bar{t}$,
    \begin{align*}
     & |Q_{t+1}^{\prime}(x,u)-Q_{t+1}^{\prime\prime}(x,u)| \\ & \qquad \leq \Pi_{s=0}^{m_t(x,u)} (1-\alpha_s) M_0 \\ & \qquad +\left[1 - \Pi_{s=0}^{m_t(x,u)} (1-\alpha_s) \right] \beta \bar{M}(1+1/\beta)/2
    \end{align*}
where $m_t(x,u):=\sum_{s=0}^t I_{\{(x_t,u_t)=(x,u)\}}$ is the number of visits to $(x,u)$ in $[0,t]$.
Since each $(x,u)$ is visited infinitely often w.p. $1$ and $\sum_s \alpha_s =\infty$, we have, for each $(x,u)$, $\Pi_{s=0}^{m_t(x,u)} (1-\alpha_s) \rightarrow 0$ as $t\rightarrow\infty$ w.p. $1$. This implies that $\bar{M}\leq \beta\bar{M}(1+1/\beta)/2<\bar{M}$ w.p. $1$, which is a contradiction.   Therefore, $M_t\rightarrow0$, w.p. $1$.

Theorem~4 in \cite{tsitsiklis1994asynchronous} shows that, for any initial condition $Q_0$, $Q_t\rightarrow\bar{Q}$, w.p. $1$. Hence,  for any  $Q_0^{\prime}\in\mathbb{Q}$, we have $|Q_t^{\prime}-\bar{Q}|_{\infty}+\sup_{Q_0^{\prime\prime}\in\mathbb{Q}}|Q_t^{\prime}-Q_t^{\prime\prime}|_{\infty}\rightarrow 0$, w.p. $1$. Therefore, $\sup_{Q_0^{\prime\prime}\in\mathbb{Q}}|Q_t^{\prime\prime}-\bar{Q}|_{\infty}\rightarrow 0$, w.p. $1$. This leads to the desired result, i.e.,
for any $\epsilon>0$ and compact $\mathbb{Q}\in\mathbb{R}^{|\mathbb{X}\times\mathbb{U}|}$, there exists $T_{\epsilon}^{\mathbb{Q}}<\infty$ such that
$$P\left[ \sup_{t \geq T_{\epsilon}^{\mathbb{Q}}} \sup_{Q_0^{\prime\prime}\in\mathbb{Q}} |Q_t^{\prime\prime}-\bar{Q}|_{\infty} \leq \epsilon  \right] \geq 1-\epsilon.$$
\end{IEEEproof}

\begin{remark}\label{rm:uc}
The Q-factors corresponding to a certain deterministic policy $\hat{\pi}$ can be learnt by modifying the recursion (\ref{eq:sQ1})-(\ref{eq:sQ2}) as follows: for $t\geq0$,
    \begin{align*}
     \hat{Q}_{t+1}(x,u)  & =   \hat{Q}_t(x,u), \qquad\mbox{for all} \ (x,u)\not=(x_t,u_t) \\
     \hat{Q}_{t+1}(x_t,u_t)  & =   \hat{Q}_t(x_t,u_t)    +  \alpha_{n_t} \left( c(x_t,u_t) \right. \\& \left.\qquad +\beta  \hat{Q}_t(x_{t+1},\hat{\pi} (x_{t+1})) - \hat{Q}_t(x_t,u_t)   \right)
    \end{align*}
where the initial condition $\hat{Q}_0$ is given and $u_t$ is chosen according to $\pi(x_t)$. Hence, the uniform convergence result in Lemma 1 also holds for the this recursion.
\end{remark}

\section{Proof of Theorem~\ref{normalQ}}\label{Qsection}

For any $\pi^{-i}\in\Delta^{-i}$, let $F_{\pi^{-i}}^i$ denote the self-mapping of $\mathbb{X}\times\mathbb{U}^i$ defined by
\begin{align*}
F_{\pi^{-i}}^i(Q^i)(x,u^i)  = &  E_{\pi^{-i}(x)} \big[ c^i\left(x,u^i,u^{-i}\right)    \\ & + \beta^i \sum_{x^{\prime}} P\left[x^{\prime}|x,u^i,u^{-i}\right] \min_{v^i} Q^i(x^{\prime},v^i) \big]
\end{align*}
for all $x,u^i$. It is well-known that $F_{\pi^{-i}}^i$ is a contraction mapping with the Lipschitz constant $\beta^i$ with respect to the maximum norm.
Recall from (\ref{eq:Qfp}) that each DM$^i$'s optimal Q-factors $Q_{\pi^{-i}}^i$  is the unique fixed point of $F_{\pi^{-i}}^i$. We also note that, during the $k-$th exploration phase, each DM$^i$ actually uses the random policy $\bar{\pi}_k^{i}$ defined as
\begin{equation}
\label{eq:pibar}
\bar{\pi}_k^j=(1-\rho^j)\pi_k^j+\rho^j\nu^j
\end{equation}
where $\nu^j$ is the random policy that assigns the uniform distribution on $\mathbb{U}^j$ to each $x$.

\begin{lemma}
\label{lm:Qappx}
For any $\epsilon>0$, there exists $T_{\epsilon}<\infty$ such that, if $T_k\geq T_{\epsilon}$, then
$$P\left[  \big|Q_{t_{k+1}}^i - Q_{\bar{\pi}_k^{-i}}^{i}\big|_{\infty} \leq \epsilon, \ \mbox{for all} \ i  \right] \geq 1-\epsilon.$$
\end{lemma}

\begin{IEEEproof}
Note that the $k-$th exploration phase starts with $x_{kT}$, which belongs to the finite state space $\mathbb{X}$, and $Q_{t_k}^i\in\mathbb{Q}^i$, where $\mathbb{Q}^i$ is compact, for all $i$. Note also that, during each exploration phase, DMs use stationary random policies of the form (\ref{eq:pibar}) and there are finitely many such joint policies. Hence, the desired result follows from Lemma~\ref{lm:ql} in Appendix~\ref{se:ap3}.\qquad
\end{IEEEproof}

\begin{lemma}
\label{lm:Qexp}
For any $\epsilon>0$, there exists $\rho_{\epsilon}>0$ such that, if  $\rho^i\leq\rho_{\epsilon}$, for all $i$, then
$$\left|Q_{\pi_k^{-i}}^i - Q_{\bar{\pi}_k^{-i}}^{i}\right|_{\infty} \leq \epsilon, \qquad \mbox{for all} \ i,  k.$$
\end{lemma}

\begin{IEEEproof}
We have
\begin{align*}
 \left| Q_{\pi_k^{-i}}^i - Q_{\bar{\pi}_k^{-i}}^{i}\right|_{\infty}
  = & \left|F_{\pi_k^{-i}}^i(Q_{\pi_k^{-i}}^i) - F_{\bar{\pi}_k^{-i}}^i(Q_{\bar{\pi}_k^{-i}}^{i}) \right|_{\infty} \\
 \leq & \left|F_{\pi_k^{-i}}^i(Q_{\pi_k^{-i}}^i) - F_{\bar{\pi}_k^{-i}}^i(Q_{\pi_k^{-i}}^{i})\right|_{\infty}  \\ & +\left|F_{\bar{\pi}_k^{-i}}^i(Q_{\pi_k^{-i}}^{i}) - F_{\bar{\pi}_k^{-i}}^i(Q_{\bar{\pi}_k^{-i}}^{i})\right|_{\infty} \\
 \leq & \left(1-\prod_{j\not=i}(1-\rho^j)\right) \times \\ & \ \left| F_{\pi_k^{-i}}^i(Q_{\pi_k^{-i}}^{i})-F_{\phi_k^{-i}}^i(Q_{\pi_k^{-i}}^{i})\right|_{\infty} \\& + \beta^i \left|Q_{\pi_k^{-i}}^i - Q_{\bar{\pi}_k^{-i}}^{i}\right|_{\infty}
\end{align*}
where $\phi_k^{-i}\in\Delta^{-i}$ is some convex combination of the policies in $\Delta^{-i}$  of the form where each DM$^j$, $j\not=i$, either uses its baseline policy $\pi_k^j\in\Pi^j$ or the uniform distribution\footnote{More precisely,  $\phi_k^{-i}=\sum_{J\subset\{1,\dots,N\}\backslash\{i\}} a_J \phi_{k,J}^{-i}$ where $a_J:=\frac{\prod_{j\in J} (1-\rho^j)\prod_{j\not\in J \cup\{i\}} \rho^j}{1-\prod_{j\not=i}(1-\rho^j)}$ and $\phi_{k,J}\in\Delta^{-i}$ is a policy such that $\phi_{k,J}^j=\pi_k^j$ for $j\in J$ and $\phi_{k,J}^j=\nu^j$ for $j\not\in J\cup\{i\}$.}. Because $(\pi_k^{-i},\phi_k^{-i})$ belongs to a finite subset of $\Pi^{-i}\times\Delta^{-i}$, an upper bound $\bar{F}<\infty$ on
$$\left| F_{\pi_k^{-i}}^i(Q_{\pi_k^{-i}}^{i})-F_{\phi_k^{-i}}^i(Q_{\pi_k^{-i}}^{i})\right|_{\infty}$$  exists, which is uniform in $(\pi_k^{-i},\phi_k^{-i})$. This results in $$\left|Q_{\pi_k^{-i}}^i - Q_{\bar{\pi}_k^{-i}}^{i}\right|_{\infty} \leq \left(1-\prod_{j\not=i}(1-\rho^j)\right) \frac{\bar{F}} {1-\beta^i}$$
which proves the lemma. \qquad
\end{IEEEproof}

Let $\bar{\delta}$ denote the minimum separation between the entries of DMs' optimal  Q-factors (with respect to the deterministic policies), defined as\footnote{To avoid trivial cases, we assume $Q_{\pi^{-i}}^i(x,v^i)\not=Q_{\pi^{-i}}^i(x,\tilde{v}^i)$ for some $i$, $x$, $v^i$, $\tilde{v}^i$, $\pi^{-i}\in\Pi^{-i}$.}
\begin{align*}
\bar{\delta}:=\min_{\begin{array}{c} \scriptstyle i,x,v^i,\tilde{v}^i,\pi^{-i}\in\Pi^{-i}: \\ \scriptstyle Q_{\pi^{-i}}^i(x,v^i)\not=Q_{\pi^{-i}}^i(x,\tilde{v}^i) \end{array}} \left|Q_{\pi^{-i}}^i(x,v^i)-Q_{\pi^{-i}}^i(x,\tilde{v}^i)\right|.
\end{align*}
We consider $\bar{\delta}$ to be an upper bound on the tolerance levels for sub-optimality, i.e., $\delta^i\in(0,\bar{\delta})$, for all $i$. In that case, we also introduce an upper bound $\bar{\rho}>0$ on the experimentation rates  such that, if $\rho^i\leq\bar{\rho}$, for all $i$, then
\begin{equation}
\label{eq:condrho}
\left|Q_{\pi_k^{-i}}^i - Q_{\bar{\pi}_k^{-i}}^{i}\right|_{\infty} < \frac{1}{2}\min\{\delta^i,\bar{\delta}-\delta^i\}, \quad\mbox{for all} \ i, \ k.
\end{equation}
Such an upper bound $\bar{\rho}>0$ exists due to Lemma~\ref{lm:Qexp}.

\begin{lemma}
\label{lm:Pi}
Suppose $\delta^i\in(0,\bar{\delta})$, $\rho^i\in(0,\bar{\rho})$, for all $i$. For any $\epsilon>0$, there exist $\bar{T}<\infty$, such that, if $T_k\geq \bar{T}$, then
$$P\left[ E_k  \right] \geq 1-\epsilon$$
where $E_k$, $k\geq0$, is the random event defined as
\begin{align*}
E_k:=\Big\{\omega\in\Omega:  \left|Q_{t_{k+1}}^i - Q_{\pi_k}^{i}\right|_{\infty} < & \frac{1}{2}\min\{\delta^i,\bar{\delta}-\delta^i\}, \\ & \quad \mbox{for all} \ i   \Big\}.
\end{align*}
\end{lemma}

\begin{IEEEproof}
The desired result follows from  Lemma~\ref{lm:Qappx} and (\ref{eq:condrho}).
\end{IEEEproof}

\subsection{Proof of part (i)}
Note that $$\omega\in E_k    \Rightarrow  \Pi_{\pi_k}=\Pi_{k+1}=\Pi_{k+1}^1\times\cdots\times\Pi_{k+1}^N.$$
Therefore, we have
\begin{equation}
\label{eq:e2e}
P\left[  \pi_{k+1} =\pi_k   |  E_k, \ \pi_k\in\Pi_{\rm eq} \right] =1, \quad\mbox{for all} \ k.
\end{equation}
Since we have a weakly acyclic game at hand, for each $\pi\in\Pi$, there exists a strict best reply path of minimum length $L_{\pi}<\infty$ starting at $\pi$ and ending at an equilibrium policy. Let $L:=\max_{\pi\in\Pi} L_{\pi}$. There exists $p_{\min}\in(0,1)$ (which depends only on $\lambda^1,\dots,\lambda^N$, and $L$) such that, for all $k$,
\begin{align}
 P\big[  \pi_{k+L} \in \Pi_{\rm eq}  \big|  E_k,\dots,E_{k+L-1},   \pi_k\not\in\Pi_{\rm eq} \big] \geq p_{\min}.\label{eq:ne2e}
\end{align}
Pick $\tilde{\epsilon}\in(0,\epsilon)$ satisfying
$$\left(\frac{(1-\tilde{\epsilon})  p_{\min}}{\tilde{\epsilon}+(1-\tilde{\epsilon})  p_{\min}}-\tilde{\epsilon}\right)(1-\tilde{\epsilon}) \geq 1-\epsilon.$$
Lemma~\ref{lm:Pi} implies the existence of $\tilde{T}<\infty$ such that, if $\min_{\ell} T_{\ell} \geq \tilde{T}$, then
\begin{equation}
\label{eq:Pi}
P\left[  E_k,\dots,E_{k+L-1}  \right] \geq 1-\tilde{\epsilon},  \quad\mbox{for all} \ k.
\end{equation}
For the rest of this part, we assume $\min_{\ell} T_{\ell}\geq \tilde{T}$. From (\ref{eq:e2e}), (\ref{eq:ne2e}), (\ref{eq:Pi}), we obtain
$$P\left[  \pi_{k+L} \in \Pi_{\rm eq}   |  \pi_k\not\in\Pi_{\rm eq}  \right] \geq p_{\min}(1-\tilde{\epsilon}), \quad\mbox{for all} \ k$$
and
\begin{equation}
\label{eq:neeqx}
P\left[  \pi_{k+L} = \cdots = \pi_k  |  \pi_k\in\Pi_{\rm eq}  \right] \geq 1-\tilde{\epsilon}, \quad\mbox{for all} \ k.
\end{equation}
This leads to the recursive inequalities
\begin{eqnarray}
p_{(n+1)L} \geq (1-\tilde{\epsilon})  [p_{nL} + p_{\min} (1-p_{nL})] \label{geometricConv}
\end{eqnarray}
where
$p_k:=P\left[  \pi_{k} \in \Pi_{\rm eq}  \right]$, for all $k$. Note that we have, for all $n$,
\begin{equation}
\label{eq:lwb}
p_{(n+1)L} -p_{nL}\geq -\tilde{\epsilon}.
\end{equation}
We rewrite (\ref{geometricConv}) as
$$
p_{(n+1)L} -p_{nL} \geq \left[\tilde{\epsilon}+(1-\tilde{\epsilon})p_{\min}\right] \left[\frac{(1-\tilde{\epsilon})  p_{\min}}{\tilde{\epsilon}+(1-\tilde{\epsilon})  p_{\min}}-p_{nL}\right].
$$
This shows that if
\begin{equation}
\label{eq:ineqxyzw}
p_{nL} \leq \frac{(1-\tilde{\epsilon})  p_{\min}}{\tilde{\epsilon}+(1-\tilde{\epsilon})  p_{\min}}-\tilde{\epsilon}
\end{equation}
we have $p_{(n+1)L} \geq p_{nL}+p_{\min}\tilde{\epsilon}$. Therefore, whenever $p_{nL}$ satisfies (\ref{eq:ineqxyzw}), it will increase by at least $p_{\min}\tilde{\epsilon}$ until it exceeds the right hand side of (\ref{eq:ineqxyzw}), which will happen in a finite number of steps. In fact, $p_{nL}$ would increase as long as $p_{nL} < \frac{(1-\tilde{\epsilon})  p_{\min}}{\tilde{\epsilon}+(1-\tilde{\epsilon})  p_{\min}}$. On the other hand, if $p_{nL} \geq \frac{(1-\tilde{\epsilon})  p_{\min}}{\tilde{\epsilon}+(1-\tilde{\epsilon})  p_{\min}}$, $p_{nL}$ cannot decrease more than $\tilde{\epsilon}$; recall (\ref{eq:lwb}).
Therefore, there exists $\tilde{n}<\infty$ such that, for all $n\geq\tilde{n}$, $$p_{nL}\geq \frac{(1-\tilde{\epsilon})  p_{\min}}{\tilde{\epsilon}+(1-\tilde{\epsilon})  p_{\min}}-\tilde{\epsilon}.$$ Finally, due to (\ref{eq:neeqx}), we have, for all $n\geq\tilde{n}$, $\ell\in\{1,\dots,L-1\}$, $$p_{nL+\ell}\geq \left(\frac{(1-\tilde{\epsilon})  p_{\min}}{\tilde{\epsilon}+(1-\tilde{\epsilon})  p_{\min}}-\tilde{\epsilon}\right)(1-\tilde{\epsilon})\geq 1-\epsilon.$$

\subsection{Proof of part (ii)} For any $\epsilon>0$, let $\tilde{T}<\infty$, $\tilde{k}<\infty$ be as in part (i). Let $\hat{k}<\infty$ be such that $\min_{k\geq\hat{k}}T_k\geq\tilde{T}$. It is straightforward to see from the proof of part (i) that, for all $k\geq\hat{k}+\tilde{k}$, we have $P\left[ \pi_{k} \in \Pi_{\rm eq} \right] \geq 1 -\epsilon.$

\subsection{Proof of part (iii)}
Pick a sequence  $\{\tilde{\epsilon}_n\}_{n\geq0}$ satisfying $\tilde{\epsilon}_n>0$, for all $n$, and
\begin{equation}
\label{eq:condet}
\sum_{n} (1-p_{\min})^{-n} \tilde{\epsilon}_n < \infty
\end{equation}
where $p_{\min}$ is as in (\ref{eq:ne2e}). Lemma~\ref{lm:Pi} implies the existence of a sequence $\{\tilde{T}_n\}_{n\geq0}$ of finite integers such that if
\begin{align} T_{nL},\dots,T_{(n+1)L-1}\geq \tilde{T}_n\label{eq:xyzz}
\end{align} then
\begin{align}
 P\left[  E_{nL},\dots,E_{(n+1)L-1}  \right]  \geq 1-\tilde{\epsilon}_n.\label{eq:xyz}
\end{align}
We assume (\ref{eq:xyzz}) (therefore (\ref{eq:xyz})) holds for all $n$. This leads to
\begin{align*}
P[  \pi_{(n+1)L} \not\in \Pi_{\rm eq} ]  \leq  (1-p_{\min}) P\left[  \pi_{nL} \not\in \Pi_{\rm eq}  \right]  + \tilde{\epsilon}_n.
\end{align*}
From this, it is straightforward to obtain
\begin{align*}
& P\left[  \pi_{(n+1)L} \notin \Pi_{e}  \right]  \\ & \qquad \leq  (1-p_{\min})^n \left( 1  + \sum_{s=0}^n (1-p_{\min})^{-s} \tilde{\epsilon}_s \right).
\end{align*}
Due to (\ref{eq:xyz}), we have, for $\ell\in\{0,\dots,L-1\}$, $$P\left[  \pi_{nL+\ell} \in \Pi_{\rm eq}  \right]  \geq (1-\tilde{\epsilon}_n) P\left[  \pi_{nL} \in \Pi_{\rm eq}  \right].$$
Therefore, for $\ell\in\{0,\dots,L-1\}$,
\begin{align*}
 & P\left[  \pi_{(n+1)L+\ell} \notin \Pi_{\rm eq}  \right] \\ & \qquad \leq (1-p_{\min})^n \left( 1  + \sum_{s=0}^n (1-p_{\min})^{-s} \tilde{\epsilon}_s \right)+\tilde{\epsilon}_{n+1}.
\end{align*}
From this and (\ref{eq:condet}), we obtain
\begin{align*}
& \sum_{k\geq1} P\left[  \pi_k \notin \Pi_{e}  \right] \\ & \ \leq L \sum_{n\geq0} \left[(1-p_{\min})^n \left( 1  + \sum_{s=0}^n (1-p_{\min})^{-s} \tilde{\epsilon}_s \right)+\tilde{\epsilon}_{n+1}\right] \\ & \ \ <\infty.
\end{align*}
Borel-Cantelli Lemma implies
\begin{align}
\label{eq:bc1}
P[\pi_k \notin \Pi_{\rm eq}, \ \mbox{for infinitely many} \ k]=0.
\end{align}
From (\ref{eq:condet}) and (\ref{eq:xyz}), we obtain $\sum_{k\geq0} P\left[ \Omega\backslash E_k\right] <\infty$. Borel-Cantelli Lemma again implies
\begin{align}
\label{eq:bc2}
P[\Omega\backslash E_k, \ \mbox{for infinitely many} \ k]=0.
\end{align}
Finally, (\ref{eq:bc1}) and (\ref{eq:bc2}) imply the desired result.

\section{Proof of Theorem~\ref{normalQ2}}\label{Qsection2}

For any $\pi=(\pi^i,\pi^{-i})\in\Pi^i\times\Delta^{-i}$,  let $F_{\pi}^i$ denote the self-mapping of $\mathbb{X}\times\mathbb{U}^i$ defined by
\begin{align*}
F_{\pi}^i(Q^i)(x,u^i)  = & E_{\pi^{-i}(x)} \big[ c^i\left(x,u^i,u^{-i}\right)   \\ &  + \beta^i \sum_{x^{\prime}} P\left[x^{\prime}|x,u^i,u^{-i}\right] Q^i(x^{\prime},\pi^i(x^{\prime})) \big]
\end{align*}
for all $x,u^i$. It is well-known that $F_{\pi}^i$ is a contraction mapping with the Lipschitz constant $\beta^i$ with respect to the maximum norm.
Let us denote the unique fixed point of $F_{\pi}^i$ by $Q_{\pi}^i$.  We also note that, during the $k-$th exploration phase, each DM$^i$ actually uses the random policy $\bar{\pi}_k^{i}$ defined as
\begin{equation}
\label{eq:pibar2}
\bar{\pi}_k^j=(1-\rho^j)\pi_k^j+\rho^j\nu^j
\end{equation}
where $\nu^j$ is the random policy that assigns the uniform distribution on $\mathbb{U}^j$ to each $x$.

\begin{lemma}
\label{lm:Qappx2}
For any $\epsilon>0$, there exists $T_{\epsilon}<\infty$ such that, if $T\geq T_{\epsilon}$, then
\begin{align*}
& P\Big[  \left|Q_{t_{k+1}}^i - Q_{(\pi_k^i,\bar{\pi}_k^{-i})}^{i}\right|_{\infty} \leq \epsilon \ \mbox{and} \\ & \qquad  \left|\hat{Q}_{t_{k+1}}^i - Q_{(\hat{\pi}_k^i,\bar{\pi}_k^{-i})}^{i}\right|_{\infty}\leq \epsilon,   \mbox{for all} \ i  \Big] \geq 1-\epsilon, \ \mbox{for all} \ k.
\end{align*}
\end{lemma}

\begin{IEEEproof}
Note that each exploration phase starts with $x_{kT}$, which belongs to a finite state space, and $Q_{kT}^i,\hat{Q}_{kT}^i\in\mathbb{Q}^i$, where $\mathbb{Q}^i$ is compact, for all $i$. Note also that, during each exploration phase, DMs use stationary random policies of the form (\ref{eq:pibar2}) and there are finitely many such joint policies. Hence, the desired result follows from Lemma~\ref{lm:ql} in Appendix~\ref{se:ap3}; see Remark~\ref{rm:uc}.\qquad
\end{IEEEproof}

\begin{lemma}
\label{lm:Qexp2}
For any $\epsilon>0$, there exists $\rho_{\epsilon}>0$ such that, if  $\rho^i\leq\rho_{\epsilon}$, for all $i$, then
\begin{align*}
& \left|Q_{(\pi_k^i,\pi_k^{-i})}^i - Q_{(\pi_k^i,\bar{\pi}^{-i}_k)}^{i}\right|_{\infty}  \leq \epsilon   \\ & \ \mbox{and}
\left|Q_{(\hat{\pi}_k^i,\pi_k^{-i})}^i - Q_{(\hat{\pi}_k^i,\bar{\pi}^{-i}_k)}^{i}\right|_{\infty} \leq \epsilon, \qquad \mbox{for all} \ i, k.
\end{align*}
\end{lemma}

\begin{IEEEproof}
We have
\begin{align*}
& \left|Q_{(\pi_k^i,\pi_k^{-i})}^i - Q_{(\pi_k^i,\bar{\pi}^{-i}_k)}^{i}\right|_{\infty}   \\
& \quad= \left|F_{(\pi_k^i,\pi_k^{-i})}^i\left(Q_{(\pi_k^i,\pi_k^{-i})}^i\right) - F_{(\pi_k^i,\bar{\pi}^{-i}_k)}^i\left(Q_{(\pi_k^i,\bar{\pi}^{-i}_k)}^{i}\right) \right|_{\infty} \\
& \quad \leq
\left|F_{(\pi_k^i,\pi_k^{-i})}^i\left(Q_{(\pi_k^i,\pi_k^{-i})}^i\right) - F_{(\pi_k^i,\bar{\pi}^{-i}_k)}^i\left(Q_{(\pi_k^i,\pi^{-i}_k)}^{i}\right) \right|_{\infty} \\
& \qquad +
\left|F_{(\pi_k^i,\bar{\pi}^{-i}_k)}^i\left(Q_{(\pi_k^i,\pi_k^{-i})}^i\right) - F_{(\pi_k^i,\bar{\pi}^{-i}_k)}^i\left(Q_{(\pi_k^i,\bar{\pi}^{-i}_k)}^{i}\right) \right|_{\infty}
\\
& \quad \leq \left(1-\prod_{j\not=i}(1-\rho^j)\right) \times \\ & \qquad \qquad\left| F_{(\pi_k^i,\pi_k^{-i})}^i\left(Q_{(\pi_k^i,\pi_k^{-i})}^i\right)-F_{(\pi_k^i,\phi_k^{-i})}^i\left(Q_{(\pi_k^i,\pi_k^{-i})}^i\right)\right|_{\infty} \\
& \qquad + \beta^i \left|Q_{(\pi_k^i,\pi_k^{-i})}^i-Q_{(\pi_k^i,\bar{\pi}^{-i}_k)}^{i} \right|_{\infty}
\end{align*}
where $\phi_k^{-i}\in\Delta^{-i}$ is some convex combination of the joint policies  of the form where each DM$^j$, $j\not=i$, either uses its baseline policy $\pi_k^j\in\Pi^j$ or the uniform distribution (as in Appendix~\ref{Qsection}). Because $(\pi_k^i,\pi_k^{-i},\phi_k^{-i})$ belongs to a finite subset of $\Pi^i\times\Pi^{-i}\times\Delta^{-i}$, an upper bound $\check{F}<\infty$ on
$$\left| F_{(\pi_k^i,\pi_k^{-i})}^i\left(Q_{(\pi_k^i,\pi_k^{-i})}^i\right)-F_{(\pi_k^i,\phi_k^{-i})}^i\left(Q_{(\pi_k^i,\pi_k^{-i})}^i\right)\right|_{\infty}$$ exists, which is uniform in $(\pi_k^i,\pi_k^{-i},\phi_k^{-i})$. This results in
$$\left|Q_{(\pi_k^i,\pi_k^{-i})}^i-Q_{(\pi_k^i,\bar{\pi}^{-i}_k)}^{i} \right|_{\infty} \leq \left(1-\prod_{j\not=i}(1-\rho^j)\right) \frac{\check{F}} {1-\beta^i}$$
which leads to the first bound. The second bound can be obtained similarly.\qquad
\end{IEEEproof}

Let $\check{\delta}$ denote the minimum separation between the entries of DMs'  Q-factors (for deterministic policies), defined as\footnote{We assume $Q_{(\pi^i,\pi^{-i})}^i(x,\pi^i(x)) \not= Q_{(\tilde{\pi}^i,\pi^{-i})}^i(x,\tilde{\pi}^i(x))$, for some $i$,  $x$, $\pi^i,\tilde{\pi}^i\in\Pi^i$, $\pi^{-i}\in\Pi^{-i}$, to avoid trivial cases.}
\begin{align*}
\check{\delta}:=\min \Big\{ &  \big|Q_{(\pi^i,\pi^{-i})}^i(x,\pi^i(x))-Q_{(\tilde{\pi}^i,\pi^{-i})}^i(x,\tilde{\pi}^i(x))\big| : \\
                            &  \quad   i,x,\pi^i,\tilde{\pi}^i\in\Pi^i, \pi^{-i}\in\Pi^{-i}, \\
                            &  \quad Q_{(\pi^i,\pi^{-i})}^i(x,\pi^i(x)) \not= Q_{(\tilde{\pi}^i,\pi^{-i})}^i(x,\tilde{\pi}^i(x))  \Big\}.
\end{align*}
We consider $\check{\delta}$ to be an upper bound on the tolerance levels for sub-optimality, i.e., $\delta^i\in(0,\check{\delta})$, for all $i$. In that case, we also introduce an upper bound $\check{\rho}>0$ on the experimentation rates  such that, if $\rho^i\leq\check{\rho}$, for all $i$, then
\begin{align}
\nonumber
 \max\Big\{ & \left|Q_{(\pi^i_k,\pi_k^{-i})}^i - Q_{(\pi^i_k,\bar{\pi}_k^{-i})}^{i}\right|_{\infty}, \\ & \left|Q_{(\hat{\pi}^i_k,\pi_k^{-i})}^i - Q_{(\hat{\pi}^i_k,\bar{\pi}_k^{-i})}^{i}\right|_{\infty} \Big\}  < \frac{1}{2}\min\{\delta^i,\check{\delta}-\delta^i\}  \label{eq:condrho2}
\end{align}
for all $i$, $k$.
Such an upper bound $\check{\rho}>0$ exists due to Lemma~\ref{lm:Qexp2}.

\begin{lemma}
\label{lm:Pi2}
Suppose $0<\delta^i<\check{\delta}$, $0<\rho^i<\check{\rho}$, for all $i$. For any $\epsilon>0$, there exist $\bar{T}<\infty$, such that, if $T_k\geq \bar{T}$, then
$$P\left[ \check{E}_k  \right] \geq 1-\epsilon$$
where $\check{E}_k$, $k\geq0$, is the random event defined as
\begin{align*}
   \check{E}_k:=\Big\{\omega\in\Omega:   \max\Big\{ & \left|Q_{t_{k+1}}^i - Q_{(\pi_k^i,\pi_k^{-i})}^{i}\right|_{\infty}, \\ & \left|\hat{Q}_{t_{k+1}}^i - Q_{(\hat{\pi}_k^i,\pi_k^{-i})}^{i}\right|_{\infty} \Big\}  \\    < \frac{1}{2} & \min\{\delta^i,\check{\delta}-\delta^i\}, \ \mbox{for all} \ i  \Big\}.
 \end{align*}
\end{lemma}

\begin{IEEEproof}
The desired result follows from Lemma~\ref{lm:Qappx2} and (\ref{eq:condrho2}).\qquad
\end{IEEEproof}

We have
\begin{equation}
\label{eq:e2e2}
P\left[  \pi_{k+1} =\pi_k   |  \check{E}_k, \ \pi_k\in\Pi_{\rm eq} \right] =1, \qquad\mbox{for all} \ k.
\end{equation}
Since we have a weakly acyclic game at hand, for each $\pi\in\Pi$, there exists a strict better reply path of minimum length $\check{L}_{\pi}<\infty$ starting at $\pi$ and ending at an equilibrium policy. Let $\check{L}:=\max_{\pi\in\Pi} \check{L}_{\pi}$. There exists $\check{p}_{\min}\in(0,1)$ (which depends only on $\lambda^1,\dots,\lambda^N$, and $L$) such that, for all $k$,
\begin{align}
 P\big[  \pi_{k+\check{L}} \in \Pi_{\rm eq}  \big|  \check{E}_k,\dots,\check{E}_{k+L-1},   \pi_k\not\in\Pi_{\rm eq} \big] \geq \check{p}_{\min}. \label{eq:ne2e2}
\end{align}
Pick $\check{\epsilon}\in(0,\epsilon)$ satisfying
$$\left(\frac{(1-\check{\epsilon})  \check{p}_{\min}}{\check{\epsilon}+(1-\check{\epsilon})  \check{p}_{\min}}-\check{\epsilon}\right)(1-\check{\epsilon}) \geq 1-\epsilon.$$
Lemma~\ref{lm:Pi2} implies the existence of $\check{T}<\infty$ such that, if $\min_{\ell}T_{\ell}\geq \check{T}$, then
\begin{equation}
\label{eq:Pi2}
P\left[  \check{E}_k,\dots,\check{E}_{k+L-1}  \right] \geq 1-\check{\epsilon}, \quad \mbox{for all} \ k.
\end{equation}
For the rest of the proof, we assume $\min_{\ell}T_{\ell}\geq \check{T}$. From (\ref{eq:e2e2}), (\ref{eq:ne2e2}), (\ref{eq:Pi2}), we obtain, for all $k$,
\begin{align*}
& P\left[  \pi_{k+\check{L}} \in \Pi_{\rm eq}   |  \pi_k\not\in\Pi_{\rm eq}  \right] \geq \check{p}_{\min}(1-\check{\epsilon}) \\
& \mbox{and} \ P\left[  \pi_{k+\check{L}} = \cdots = \pi_k  |  \pi_k\in\Pi_{\rm eq}  \right] \geq 1-\check{\epsilon}.
\end{align*}
This leads to the recursive inequalities
\begin{eqnarray}
p_{(n+1)\check{L}} \geq (1-\check{\epsilon})  [p_{n\check{L}} + \check{p}_{\min} (1-p_{n\check{L}})], \quad n\geq0 \label{geometricConv2}
\end{eqnarray}
where
$p_k:=P\left[  \pi_{k} \in \Pi_{\rm eq}  \right]$. Note that these inequalities are similar to (\ref{geometricConv}) and by similar reasoning, there exists $\check{n}<\infty$ such that, for all $n\geq\check{n}$ and $\ell\in\{1,\dots,L-1\}$, $$p_{n\check{L}+\ell}\geq \left(\frac{(1-\check{\epsilon})  \check{p}_{\min}}{\check{\epsilon}+(1-\check{\epsilon})  \check{p}_{\min}}-\check{\epsilon}\right)(1-\check{\epsilon})\geq 1-\epsilon.$$
This proves part (i).
The proofs of part (ii)-(iii) are analogous to the proofs of part (ii)-(iii) of Theorem~\ref{normalQ}, respectively.

\end{document}